\documentclass[12pt]{amsart}

\usepackage{amsmath}
\usepackage{amsthm}
\usepackage{amssymb}

\usepackage{pstricks}

\usepackage{ifthen}

\input diagrams

\newcommand{\trt}{2}
\newcommand{\lbl}[1]
{\ifthenelse{\trt=2}{\quad\label{#1}\ref{#1}}{\label{#1}}}

\newfont{\rmm}{cmr10 scaled 1000}
\newfont{\itt}{cmsl10 scaled 1000}
\newfont{\nazad}{cmff10 scaled 1000}

\newcommand{\bff}{\bf}

\newcommand{\rref}[1]{(\ref{#1})}

\theoremstyle{plain}

\newtheorem{theo}{Theorem}[section]
\newtheorem{lemm}[theo]{Lemma}
\newtheorem{prop}[theo]{Proposition}
\newtheorem{coro}[theo]{Corollary}

\theoremstyle{definition}

\newtheorem{defi}[theo]{Definition}
\newtheorem{rema}[theo]{Remark}

\theoremstyle{plain}
\newtheorem*{A}{Theorem A}
\newtheorem*{B}{Theorem B}

\def\novidok{S.P.Novikov,
\emph{Mnogoznachnye funktsii i funktsionaly.
Analog teorii Morsa}, Doklady AN SSSR
{\bf 260}
(1981),
 31-35.
   English translation:
 S.P.Novikov. \emph{ Many-valued functions
 and functionals. An analogue of Morse theory},
 Sov.Math.Dokl.
{\bff 24}
(1981),
       222-226. }

\newcommand{\armaz}
{M.Artin, B.Mazur,
\emph{On periodic points},
Annals of Math. {\bff 102} (1965), 82--99.
}

\newcommand{\chapman}
{ T.A.Chapman,
\emph{ Topological invariance of Whitehead torsion},
American J. of Math.
{\bff 96}
(1974),
    488 - 497
}

\newcommand{\bhs}
{H.Bass, A.Heller, R.G.Swan,
\emph{The Whitehead group of a polynomial extension},
Inst. Hautes Etudes Sci. Publ. Math. {\bff 22} (1964),
61--79
}

\newcommand{\farran}{ M. Farber and A. A. Ranicki,
\emph{ The Morse-Novikov theory of circle-valued functions
and noncommutative localization,}
 e-print dg-ga/9812122,  Proc. 1998 Moscow Conference 
for the 60th Birthday of S. P. Novikov, tr.
 Mat. Inst. Steklova, {\bf 225}, 1999,
381 -- 388}

\newcommand{\fried}{  D.Fried,
\emph{Homological Identities for closed orbits},
Inv. Math. {\bff 71}, (1983) 419--442.
}

\newcommand{\friedtwi}{  D.Fried,
\emph{Periodic points and twisted coefficients},
Lect. Notes in Math.,
{\bff 1007},
(1983)
261--293.
}

\newcommand{\friednewzeta}{  D.Fried,
\emph{Flow equivalence, hyperbolic systems and a new zeta function
for flows},
Comm. Math. Helv.,
{\bff 57},
(1982)
237--259.
}

\newcommand{\fuller}{  F.B.Fuller,
\emph{An index of fixed point type for periodic orbits},
Amer. J.Math
{\bff 89},
(1967)
133--148
}

\def\gnXCIV
{  R.Geoghegan, A.Nicas,
\emph{Trace and torsion in the theory of flows},
Topology,
{\bff 33},
(1994)
683--719
}

\newcommand{\hulee}
{M.Hutchings, Y-J.Lee
\emph{              Circle-valued Morse theory, Reidemeister torsion
and Seiberg-Witten invariants of 3-manifolds},
Topology,
{\bff 38},
(1999),
861 -- 888,
e-print dg-ga/9612004 3  Dec 1996.}

\newcommand{\milncyccov}
{J.Milnor,
\emph{ Infinite cyclic coverings},
In: Conference on the topology of manifolds,
(1968)}

\newcommand{\milnWT}
{J.Milnor,
\emph{ Whitehead Torsion},
Bull. Amer. Math. Soc.
{\bff 72}
(1966),
358 - 426.
}

\def\ranXCIX
{ A.Ranicki,
\emph{The algebraic construction 
of the Novikov complex of a circle-valued Morse function},
e-print math.AT/9903090, to appear in "Mathematische Annalen"}

\def\schuetzC
{  D.Sch\"utz,
\emph{Gradient flows of closed 1-forms and their closed orbits},
math.DG/0009055, to appear in  "Forum Mathematicum"}

\def\schuetzCI
{  D.Sch\"utz,
\emph{One-parameter fixed point theory and gradient flows
of closed 1-forms},
math.DG/0104245}

\newcommand{\smdyn}
{  S.~Smale,
\emph{Differential dynamical systems},
 Bull. Amer. Math. Soc. {\bff 73} (1967)
747--817.}

\newcommand{\patou}
{ A.V.Pajitnov, \emph{ On the Novikov complex for rational Morse forms},
preprint:
Institut for Matematik og datalogi, Odense Universitet 
Preprints 1991, No 12, Oct. 1991; \quad
Annales de la Facult\'e de Sciences de 
Toulouse {\bff 4}  (1995), 297--338.
}

\newcommand{
\pastpet}
{  A.V.Pajitnov,
\emph{   
Ratsional'nost' granichnyh operatorov v komplekse 
Novikova v sluchae obschego polozheniya},
e-print:
 dg-ga/9603006 14 Mar 96,\quad
Algebra i Analiz
{\bff 9}, no.5 (1997),  92--139.
English translation:
{\it   The incidence coefficients in the Novikov complex are
generically rational functions},
Sankt-Petersbourg Mathematical Journal
\textbf{9}
(1998),
no. 5, p. 969 -- 1006,
}

\newcommand{\pafest}
{  A.V.Pajitnov,
\emph{   Prostoi gomotopicheskii tip 
kompleksa Novikova i $\zeta$-funktsiya Lefschetza
gradientnogo potoka}, 
e-print
dg-ga/970614 9 July 1997,
,
Uspekhi Mat. Nauk
\textbf{54}
(1999)
no. 1,     117 -- 170,
English translation:
{\it Simple homotopy type of Novikov complex
and $\zeta$-function of the gradient flow},
Russian Mathematical Surveys,  1999}

\newcommand{\paadv}
{  A.V.Pajitnov,
\emph{   $C^0$-generic properties of boundary operators in the Novikov
complex },
e-print: math.DG/9812157, 29 Dec 1998,
Advances in Mathematical Sciences, vol. 197, 1999, p.29 -- 117;
}

\newcommand{\pawitt}
{  A.V.Pajitnov,
\emph{    Closed orbits of gradient 
flows and logarithms of non-abelian Witt vectors},
e-print math.DG/9908010, 2 Aug. 1999
  K-theory, Vol. 21 No. 4, 2000}

\newcommand{\atiyahandmacdo}
{   M.F.Atiyah, I.G.Macdonald
\emph{Introduction to commutative algebra
},
Addison-Wesley,   1969.
}

\newcommand{\dold}
{   A.Dold
\emph{Lectures on Algebraic Topology
},
Springer,  1972.
}

\newcommand{\franks}
{   J.Franks
\emph{Homology and dynamical systems},
CBMS Reg. Conf. vol. 49, AMS, Providence 1982.
}

\newcommand{\milnhcob}
{   J.~Milnor,
\emph{Lectures on the
h-cobordism theorem},
 Princeton University
Press,
 1965.
}

\newcommand{\munkres}
{  J.R.Munkres,
\emph{Elementary differential toplogy},
Annals of Math. Studies,
Vol.54, Pinceton
 1963.}

\newcommand{\massey}
{   W.Massey,
\emph{ Homology and cohomology theory
},
 Marcel Dekker, 1978.
}
\newcommand{\milnmt}
{   J.~Milnor,
\emph{ Morse theory
},
 Princeton University Press, 1963.
}

\newcommand{\spanier}
{ E.H.Spanier,
\emph{Algebraic topology},
  McGraw-Hill, 
( 1966)}

\renewcommand{\a}{\alpha}
\renewcommand{\b}{\beta}
\newcommand{\g}{\gamma}
\renewcommand{\d}{\delta}
\newcommand{\e}{\epsilon}
\newcommand{\ve}{\varepsilon}
\newcommand{\z}{\zeta}

\renewcommand{\l}{\lambda}

\newcommand{\s}{\sigma}

\newcommand{\G}{\Gamma}
\newcommand{\D}{\Delta}

\renewcommand{\L}{\Lambda}

\renewcommand{\AA}{{\mathcal A}}
\newcommand{\BB}{{\mathcal B}}
\newcommand{\CC}{{\mathcal C}}
\newcommand{\DD}{{\mathcal D}}
\newcommand{\EE}{{\mathcal E}}
\newcommand{\FF}{{\mathcal F}}
\newcommand{\GG}{{\mathcal G}}
\newcommand{\HH}{{\mathcal H}}

\newcommand{\JJ}{{\mathcal J}}

\newcommand{\LL}{{\mathcal L}}
\newcommand{\MM}{{\mathcal M}}
\newcommand{\NN}{{\mathcal N}}

\newcommand{\RR}{{\mathcal R}}

\newcommand{\TT}{{\mathcal T}}

\newcommand{\ZZ}{{\mathcal Z}}

\newcommand{\NNN}{{\mathbf{N}}}

\newcommand{\QQQ}{{\mathbf{Q}}}
\newcommand{\RRR}{{\mathbf{R}}}

\newcommand{\ZZZ}{{\mathbf{Z}}}

\newcommand{\gC}{{\mathfrak{C}}}

\newcommand{\id}{\text{id}}

\newcommand{\Wh}{\text{\rm Wh }}
\newcommand{\Ker}{\text{\rm Ker }}

\newcommand{\ind}{\text{\rm ind}}

\renewcommand{\Im}{\text{\rm Im }}
\newcommand{\supp}{\text{\rm supp }}
\newcommand{\Int}{\text{\rm Int }}

\newcommand{\Fix}{\text{\rm Fix }}

\newcommand{\bere}{\begin{rema}}
\newcommand{\bede}{\begin{defi}}

\renewcommand{\beth}{\begin{theo}}
\newcommand{\bele}{\begin{lemm}}
\newcommand{\bepr}{\begin{prop}}
\newcommand{\beeq}{\begin{equation}}
\newcommand{\bega}{\begin{gather}}
\newcommand{\been}{\begin{enumerate}}

\newcommand{\beco}{\begin{coro}}

\newcommand{\beal}{\begin{aligned}}

\newcommand{\enre}{\end{rema}}

\newcommand{\enco}{\end{coro}}
\newcommand{\enpr}{\end{prop}}
\newcommand{\enth}{\end{theo}}
\newcommand{\enle}{\end{lemm}}
\newcommand{\enen}{\end{enumerate}}
\newcommand{\enga}{\end{gather}}
\newcommand{\eneq}{\end{equation}}
\newcommand{\enal}{\end{aligned}}

\newcommand{\bq}{\begin{equation}}
\newcommand{\bqq}{\begin{equation*}}

\renewcommand{\leq}{\leqslant}
\renewcommand{\geq}{\geqslant}

\newcommand{\vide}{  \emptyset  }
\newcommand{\bu}{\bullet}

\newcommand{\mx}{\mbox}

\newcommand{\wi}{\widetilde}

\newcommand{\ove}{\overline}

\newcommand{\wh}{\widehat}

\newcommand{\sm}{\setminus}

\newcommand{\sbs}{\subset}

\newcommand{\subs}{\subsection}

\newcommand{\lb}{\label}

\newcommand{\st}[1]{\overset{\rightsquigarrow}{#1}}

\newcommand{\stexp}[1]{{#1}^{\rightsquigarrow}}

\newcommand{\stv}{\stexp {(-v)}}
\newcommand{\stu}{\stexp {(-u)}}
\newcommand{\stw}{\stexp {(-w)}}

\newcommand{\vflesh}{v\!\da}

\newcommand{\da}{\downarrow}

\newcommand{\tens}[1]{\underset{#1}{\otimes}}

\newcommand{\dow}{\pr_0 W}

\newcommand{\daw}{\pr_1 W}

\newcommand{\dwmokm}{\daw^{(k-1)}}
\newcommand{\dwmokmm}{\daw^{(k-2)}}

\newcommand{\wmok}{\Wmok}
\newcommand{\Wmok}{W^{\langle k\rangle}}

\newcommand{\wmokm}{W^{\langle k-1\rangle}}

\newcommand{\ata}{almost~ transversality~ assumption}

\newcommand{\sut}{~such~that~}

\newcommand{\wrt}{~with respect to}
\newcommand{\ho}{homomorphism}

\newcommand{\ma}{manifold}
\newcommand{\nei}{neighborhood}
\newcommand{\dfm}{diffeomorphism}

\newcommand{
\sma}{submanifold}

\newcommand{\noconf}{~there~is~no~possibility~of~confusion}

\newcommand{\TA}{transversality condition}

\newcommand{\heq}{homotopy equivalence}

\newcommand{\heeq}{homology equivalence}

\newcommand{\glvf}{gradient-like vector field}

\newcommand{\su}{subsection}

\newcommand{\Prf}{{\it Proof.\quad}}

\newcommand{\smo}{C^{\infty}}

\newcommand{\chart}{\Phi_p:U_p\to B^n(0,r_p)}
\newcommand{\atlas}{\{\Phi_p:U_p\to B^n(0,r_p)\}_{p\in S(f)}}

\newcommand{\fcob}{f:W\to[a,b]}

\newcommand{\indl}[1]{{\scriptstyle{\text{\rm ind}\leqslant {#1}~}}}
\newcommand{\inde}[1]{{\scriptstyle{\text{\rm ind}      =   {#1}~}}}
\newcommand{\indg}[1]{{\scriptstyle{\text{\rm ind}\geqslant {#1}~}}}

\newcommand{\pr}{\partial}

\newcommand{\qt}{\hfill\triangle}
\newcommand{\qs}{\hfill\square}

\newcommand{\pa}{\vskip0.1in}

\renewcommand{\(}{\big(}
\renewcommand{\)}{\big)}

\newcommand{\liminv}{\underset {\leftarrow}{\lim}}

\newcommand{\kpr}{K_r^+}

\newcommand{\kmr}{K_r^-}

\newcommand{\kpd}{K_r^+(\d)}

\newcommand{\kmd}{K_r^-(\d)}

\begin{document}

\title[Counting closed orbits]
{Counting closed orbits of gradients of circle-valued maps}
\author{A.V.Pajitnov}
\address{UMR 6629 CNRS, Universit\'e de Nantes\\
D\'epartement de Math\'ematiques\\
2, rue de la Houssini\`ere, 44072, Nantes Cedex France}
\email{pajitnov@math.univ-nantes.fr}
\begin{abstract}
Let $M$ be a closed connected manifold,
$f:M\to S^1$
 be a Morse map, belonging to an indivisible integral 
class
$\xi\in H^1(M)$,
$v$ be an $f$-gradient
satisfying the transversality condition.
The Novikov construction associates to these
data a chain complex $C_*=C_*(f,v)$. 
The first main 
result of the paper is the construction
of a {\it functorial} chain homotopy equivalence from
$C_*$ to the
completed simplicial chain complex of the 
infinite cyclic covering of $M$,
corresponding to $\xi$.
The second main result
states that the torsion of this chain homotopy 
equivalence is equal to the Lefschetz zeta function 
of the gradient flow for any gradient-like
 vector field $v$ satisfying the transversality condition
and having only hyperbolic closed orbits.
\end{abstract}
\maketitle

\section{Introduction}
\label{s:intro}

The subject of the present paper belongs to  the intersection
of two domains of
topology:
Morse-Novikov theory of circle-valued functions, and the theory of
dynamical
zeta functions.
We start by a quick recollection of the basics
of the Morse-Novikov theory, and
continue with the background up to Section \ref{su:stat_res}
which contains the statement of the results of the paper.

\subsection{Morse complex}
\label{su:morse_complex}

One of the basic constructions in the classical Morse theory
is that of {\it Morse complex}.
Starting with a Morse function
$g:M\to\RRR$
on a closed \ma~$M$, and a
 gradient-like vector field $v$ for $g$,
this construction
gives a chain complex $C_*(g,v)$,
which computes
the homology of $M$
(here the vector field $v$ must satisfy the \TA,
that is the stable and unstable manifolds of all critical points
intersect transversally).
The group $C_k(g,v)$
is a free abelian group, freely generated by
the critical points of $g$
of index $k$,
and the boundary operator is defined via counting flow lines of $v$
joining critical points of $g$.

\subsection{Novikov complex}
\label{su:novikov_complex}

In the beginning of 80s this construction was generalized
by S.P.Novikov \cite{novidok}
to the case of circle-valued Morse functions.
The input data for this construction is a circle-valued Morse function
$f:M\to S^1$
on a closed connected \ma~$M$
and a gradient-like vector field $v$ for $f$. 
As before, we  assume that $v$ satisfies  the \TA.
The result is a chain complex
$C_*(f,v)$
of free modules over a ring
$\wh L=\ZZZ[[t]][t^{-1}]$
of Laurent series in one variable with integral coefficients
and finite negative part.
The construction of the Novikov complex 
is recalled in  Section \ref{su:nov_basics} in   details,
here we just mention that $C_k(f,v)$
is freely generated over $\wh L$ by the critical points of $f$.
The homology of $C_*(f,v)$
is described by the following isomorphism:
\bq\label{f:first_iso}
H_*(C_*(f,v))\approx H_*(\bar M)\tens{L}\wh L
\end{equation}
Here $\bar M$
is the infinite cyclic cover of $\bar M$
induced by $f$ from $\RRR\to S^1$,
and $L=\ZZZ[t,t^{-1}]$
(we assume that the homotopy class of $f$ is indivisible
in $H_1(M,\ZZZ)$).
The proof of the isomorphism (\ref{f:first_iso})
is contained in \cite{patou}.
This paper contains actually some more than the proof of the above
isomorphism.
Namely, we gave there an explicit construction of a chain homotopy
equivalence
\bq
\lb{f:first_h_eq}
\phi:
C_*(f,v)\to C_*^\D(\bar M)\tens{L}\wh L
\end{equation}
(here $C^\D(\bar M)$
stands for the simplicial chain complex of $\bar M$).

\subsection{Dynamical  zeta functions}
\label{su:zeta_functions}

Now we turn to dynamical systems.
Let us start with a continuous map
$h:X\to X$.
In order to investigate the sets
of periodic points of $h$ of period $n$
for $n\to\infty$
Artin and Mazur
\cite{armaz}
introduced a zeta function
which encodes the information about all the periodic
points into a single power series in one variable.
This zeta function and its generalizations
have been intensively studied.
The Artin-Mazur zeta function
has a homotopy counterpart --
the {\it Lefschetz zeta function }
introduced by S.Smale \cite{smdyn}.
Here is the definition.
Let $\Fix h$
denote the set of
fixed points of $h$.
Assume  that for any $n$
the set
$\Fix h^n$
 is finite.
Set
\bq
L_k(h)=\sum_{a\in \Fix h^k} \nu(a)
\end{equation}
where $\nu(a)\in\ZZZ$
is the index of the fixed point $a$.
Define the Lefschetz zeta function by the following formula:

\begin{equation}\label{f:lefzet}
\zeta_{L}(t)=\exp \bigg(\sum_{k=1}^\infty \frac {L_k(h)}k t^k\bigg)
\end{equation}
Here is the formula for $\zeta_L$
in terms of the homology invariants of  $h$:
\begin{equation}\label{f:detzet}
\zeta_{L}(t)= \prod_i \det(I-th_{i})^{(-1)^{i+1}}
\end{equation}
where
$h_{i}$
stands for
the \ho~induced in $H_i(X)$ by $h$.

Proceeding to the dynamical systems generated by
flows on manifolds, one would expect that the theory
described above generalizes to this setting.
For every flow on a smooth manifold there should be  some power
series of the type (\ref{f:lefzet}), computable in homotopical terms.
Such results were obtained in  80s in many important cases
(see  \cite{fried})
but only for non-singular flows.
We shall cite here only the simplest
result in this direction, which goes back to
\cite{milncyccov}.
We shall reformulate it in terms which we use 
later, so we start by introducing the necessary definitions.

\subsection{Eta function and zeta function}
\label{su:eta_zeta}

Let 
$M$
be a closed connected \ma, and
$f:M\to S^1$
be a Morse map  \sut~$[f]\in H_1(M,\ZZZ)$
is indivisible.
Let
$v$ be a \glvf~for $f$.
We shall study the dynamics of the flow generated by $(-v)$.
As before, we assume that $v$ satisfies \TA. Assume moreover, that 
all closed orbits of $(-v)$ are hyperbolic. 
A natural numerical  invariant associated to a
non-degenerate closed orbit $\g$ of $(-v)$ is the
{\it Fuller index} (\cite{fuller})
$\iota_F(\g)=\frac {\ve(\g)}{m(\g)}$
where $\ve(\g)\in\{1, -1\}$
is the Poincar\'e index of $\g$, and $m(\g)\in\NNN $
is the multiplicity of $\g$.
For each closed orbit $\g$ of $(-v)$
define $n(\g)$
to be the integer, opposite to the
winding number of $f\circ \g$.
In another words, if
$[\g]\in H_1(M)$
stands for the homology class of $\g$, then
$f(\g)=-n(\g)\cdot {\bf 1}$,
where ${\bf 1}$
is the generator of $H_1(S^1)$
corresponding to the positive orientation.
Introduce the {\it Lefschetz eta function} by the following formula:
\bq\label{f:zeta_def}
\eta_L(-v)
=
\sum_{\g\in Cl(-v)} \iota_F(\g) t^{n(\g)}
\end{equation}
where
$Cl(-v)$ denotes the set of all the closed orbits of $(-v)$
(we identify closed orbits which are obtained one from another by reparametrization).
Using our  conditions on $v$ 
it is not difficult to check that the
expression in the right hand side of (\ref{f:zeta_def})
is a well-defined  power series in $t$ with rational coefficients and
vanishing constant term.
Therefore  the power series
\bq
\zeta_L(-v)= \exp(\eta_L(-v))
\end{equation}
is well defined and is again in  $\wh L_\QQQ$.
This power series is called {\it Lefschetz zeta function (of the flow
generated by  $(-v)$)}.
\bele
The power series $\zeta_L(-v)$
is in $\wh L$.
\enle
\Prf
A well-known computation
(see \cite{franks}, Prop. 5.19)
gives
\bq
\zeta_L(-v)
=
\prod_{\g\in ClPr(-v)} \big(1+\varepsilon_1(\g)t^{n(\g)}\big)^{\varepsilon_2(\g)} 
\label{f:fac_zeta}
\end{equation}
where $ClPr(-v)$
stands for the set of {\it prime} closed orbits, and
$\varepsilon_1(\g), \varepsilon_2(\g)\in\{1,-1\}$.
The righthand side of (\ref{f:fac_zeta})
is obviously in $\wh L$. $\qs$

\subsection{Counting closed orbits of a non-singular map $M\to S^1$}
\label{su:non_sing}

Assume now that  $f$ has no critical points, so that
$f$ is a fibration over $S^1$.
In this case
$H_*(\bar M)\tens{L}\wh L$
vanishes
and we can associate to $M$
a (suitably defined) Whitehead torsion.
We shall work with the following versions of $K_1$-groups:
\begin{align*}
\ove{K}_1(\wh L)
&
=
K_1(\wh L)/\{0, [-1]\},
\\
\Wh(\wh L)
&
=
K_1(\wh L)/T,
~\mbox{~where ~}
T=\{\pm t^n~|~n\in\ZZZ\}.
\end{align*}
(We denote by $\wh L^\bu$ the multiplicative group
of all units in $\wh L$, and by $[x]$ the image in
$K_1(\wh L) $
of the element 
$x\in \wh L^\bu$.)
The ring $\wh L$
is euclidean, and therefore
the determinant map
defines an isomorphism
$\det: K_1(\wh L)  \to \wh L^\bu$
and also an isomorphism
$$
\Wh(\wh L)=K_1(\wh L)/T\rTo^{\det} W
$$
where
$$
W=\{1+\sum_{i>0}a_it^i~|~a_i\in\ZZZ\}
$$
is the multiplicative group of power series with the constant term $1$
(such power series is called {\it Witt vector}).
Let
$\D$
be any $C^1$ triangulation of $M$, then
it lifts to a $\ZZZ$-invariant 
triangulation 
of $\bar M$, and the simplicial chain complex
$C_*^\D(\bar M)$
obtains a natural base, 
well defined up to the action of the elements
$\pm t^n, n\in\ZZZ$.
Therefore the torsion of the acyclic complex
$C_*^\D(\bar M)$
is well defined as an element of
$\Wh(\wh L)$.
It will be denoted by $\tau(M)$.
The theorem proved in \cite{milncyccov}
says that
\bq
\label{f:tors_miln}
\det(\tau(M))
=
(\zeta_L(-v))^{-1}
\end{equation}
Note that the left hand side 
depends only on $M$ and $[f]\in H^1(M,\ZZZ)$
but not on the particular choice of $v$.
Generalizations of this formula to the case when $f:M\to S^1$
has critical points were obtained only very recently
in \cite{hulee}, \cite{pafest}.
These two papers consider various particular cases of the 
problem; we shall explain the results in the next two subsections.

\subsection{The acyclic case: the Hutchings-Lee formula}
\label{su:hut_lee}

The first generalization
of the formula (\ref{f:tors_miln})
to the case of Morse maps with critical points
was obtained by M.Hutchings and Y-J.Lee \cite{hulee}.
Their formula for  $\zeta_L(-v)$
contains an additional term depending on the
Novikov complex.
Let $\GG(f)$
denote the set of gradient-like vector fields $v$, \sut~that
the \TA~holds and every closed orbit of $v$ is hyperbolic
(such gradient-like vector fields 
will be called {\it Kupka-Smale gradients}).
Set 
$$
L_\QQQ=\QQQ[t,t^{-1}],\quad
\wh L_\QQQ=\QQQ[[t]][t^{-1}]
$$
Assume that
$C_*(f,v)\tens{L}\wh L_\QQQ$
 and $H_*(\bar M)\tens{L}\wh L_\QQQ$
are acyclic, and denote by
$\tau_{Nov}\in \ove K_1(\wh L_\QQQ)/T $, resp.
by $\tau_M\in \ove K_1(\wh L_\QQQ)/T $
the torsions of these complexes.
Then (\cite{hulee}, Theorem 1.12) for every $v\in \GG(f)$
\bq\label{f:th_hulee}
\det(\tau_{Nov}/
\tau_M  )    
=
\zeta_L(-v)
\end{equation}

\subsection{The non-acyclic case}
\label{su:gen_c}

A formula for $\zeta_L(-v)$
without any acyclicity assumption was obtained in \cite{pafest}.
The methods of this work are different from these of
\cite{hulee}; they were introduced
earlier \cite{pastpet} to prove
the generic rationality of the 
boundary operators in the
Novikov complex.
In \cite{pafest} I proved that there is a chain homotopy equivalence
\bq\label{f:zeta_pafest}
\psi
: C_*(f,v)\to C_*^\D(\bar M)\tens{L}\wh L
\end{equation}
~\sut
~
\bq\label{f:tors_p}
\det\tau(\psi) = (\zeta_L(-v))^{-1}.
\end{equation}
However the class of gradients $v$
for which the formula (\ref{f:zeta_pafest})
was proved  is  smaller than the class of all Kupka-Smale gradients.
Namely the formula 
\rref{f:zeta_pafest}
holds for every \glvf~$v$
from a $C^0$-open-and-dense subset of $\GG(f)$.

It is natural to ask if the formula
\rref{f:tors_p}
is true for {\it every }
$f$-gradient in $\GG(f)$.
This question was posed in my paper \cite{pawitt}
(in a more general non-abelian setting)
and was the starting point for the present work.
The theorem B below (\pageref{t:B})
provides the positive answer, at least for the abelian case.

The other question suggested by the formula \rref{f:tors_p}
is whether the chain homotopy equivalence $\psi$
is homotopic to the chain equivalence
\rref{f:first_h_eq}
(the constructions of the two equivalences are different 
from one another).
This is related to the question, 
whether there is a {\it canonical}
homotopy equivalence
$C_*(f,v)\to C_*^\D(\bar M)\tens{L}\wh L$,
\sut~the formula
(\ref{f:tors_p})
holds. This question was 
posed to me 
in 1998 by M.Kervaire and V.Turaev.
Their question
 is answered by Theorem A
below, and  this  precise form of the answer 
was conjectured by V.Turaev.
The methods developed in the present paper for the proof
of the functoriality of 
the equivalence \rref{f:first_h_eq}
lead also to a quick proof that $\phi$ and $\psi$
are homotopic.

\subsection{Statement of results}
\label{su:stat_res}

The main aim of the present paper
is to obtain an analog of the formula
(\ref{f:tors_p})
for {\it any} $v\in \GG(f)$. 
The formula which we obtain (Theorem B)
generalizes both the results of \cite{hulee} and \cite{pafest},
and provides the answer in the general case.
We start by studying in details the chain homotopy equivalences
between $C_*(f,v)$
and 
$C_*^\D(\bar M)\tens{L}\wh L$.
It  is  clear from the preceding that these chain homotopy
equivalences are important  geometric objects, 
related to both the Novikov complex
and the dynamics of the gradient flow.
In the first part  of the present paper
we prove that the construction
of the  chain homotopy equivalence
$C_*(f,v)\to
C^\D_*(\bar M)\tens{L}\wh L$
given in \cite{patou}
is  {\it functorial}
(see Theorem A).
In the second part
(Theorem B)
we prove that the torsion of this chain homotopy equivalence
equals to the Lefschetz zeta function of the gradient flow for
any $v\in \GG(f)$.
In the course of the proof we show also that the chain
homotopy equivalences
$C_*(f,v)\to
C^\D_*(\bar M)\tens{L}\wh L
$
constructed in \cite{patou} and \cite{pafest}
are chain homotopic
(this is done in Subsection \ref{su:c0_gener}).
Now we proceed  to precise statements.
The next definition is convenient to state the
functoriality property.

\bede\label{d:m_triple}
{\it A Morse-Novikov triple}
(or {\it $MN$-triple } for short)
is a triple
$(M,f,v)$
where $M$ is  a closed connected manifold,
$f:M\to S^1$
a Morse map, \sut~$\xi(f)=f_*:H_1(M)\to H_1(S^1)\approx \ZZZ$
is indivisible, and
$v$  a \glvf,
satisfying \TA.         \end{defi}

Let $(M,f,v)$ be an $MN$-triple.
The Novikov construction, which is standard by now
(we recall it in Section \ref{su:nov_basics})
associates to each Morse-Novikov triple a chain complex
$C_*(f,v)$
of $\wh L$-modules, \sut~
$C_k(f,v)$
is freely generated over $\wh L$ by the set $S_k(f)$
of the critical points of $f$ of index $k$.

Let
$
(M_1,f_1,v_1),\
(M_2,f_2,v_2)
$
be two $MN$-triples.
Denote
$\xi(f_1)$
by $\xi_1$,
and
$\xi(f_2)$
by $\xi_2$.
Let $\bar M_1, \bar M_2$
be the infinite cyclic covers
corresponding to $\xi_1$, resp. $\xi_2$.
Let
$g:M_1\to M_2$
be a diffeomorphism, satisfying the following condition:
\bq\label{f:cond_func}
\a)~g_*(v_1)=v_2; \quad \b)~g^*(\xi_2)=\xi_1
\end{equation}
Let 
$\bar g: \bar M_1\to \bar M_2$
be a lift of $g$. Then $\bar g$
is a diffeomorphism 
which commutes with the action of $\ZZZ$,
and we obtain an $\wh L$-isomorphism
\bqq
\bar g_{\sharp}
:
C_*^s(\bar M_1)\tens{ L}\wh L
\to
C_*^s(\bar M_2)\tens{ L}\wh L
\end{equation*}
(here and elsewhere we denote by $C^s(X)$
the singular chain complex of $X$).
Also  $\bar g$ induces an isomorphism
$
\bar g_{!}
:
C_*(f_1,v_1)
\to
C_*(f_2,v_2)
$
of Novikov complexes
(see \ref{su:functo}).

\begin{A}
To each Morse-Novikov  triple
$(M,f,v)$
is  associated a chain homotopy equivalence
\bq\label{f:Phi}
\Phi=\Phi(M,f,v): C_*(f,v)\to C_*^s(\bar M)\tens{ L}\wh L
\end{equation}
which is functorial in the following sense:

For any two Morse-Novikov  triples
$
(M_1,f_1,v_1),\
(M_2,f_2,v_2)
$
and
 a diffeomorphism $g:M_1\to M_2$
satisfying the condition
(\ref{f:cond_func})
the following diagram is chain homotopy commutative:

\begin{diagram}[LaTeXeqno]
C_*(f_1,v_1) & \rTo^{\bar g_{!}} & C_*(f_2,v_2)
\\
\dTo_{\Phi_1} & & \dTo_{\Phi_2}
\quad\label{f:squa_theo}  \\
C_*^s(\bar M_1)\tens{ L}\wh L 
&
 \rTo^{\bar g_\sharp } & C_*^s(\bar
M_2)\tens{L}\wh L 
\end{diagram}
(here $\Phi_i=\Phi(M_i,f_i,v_i)$ for $i=1,2$, and
$\bar g:\bar M_1\to \bar M_2$
is any lift of $ g$).
\end{A}

\bere
The diffeomorphism
$\bar g$
(and the maps
$\bar g_!, \bar g_\sharp$)
are defined by $g$ uniquely up to multiplication by $t^n, n\in\ZZZ$.
\enre

If we replace  the complex of singular chains $C_*^s(\bar M)$
in 
Theorem A 
by the simplicial chain complex $C_*^\D(\bar M)$
it will be possible to consider the 
torsion of the resulting chain \heq.
Let us  recall the corresponding notions.
Let $X$ be a topological space.
Let $\pi:\wi X\to X$
be a regular covering with structure group $H$.
The singular chain complex
$C_*^s(\wi X)$ is free over $\ZZZ H$.
If $X$ is a simplicial complex
then $\wi X$
inherits a natural $H$-invariant
triangulation from $X$,
and the simplicial chain complex
$C_*^\D(\wi X)$
is a chain complex of free
finitely generated
$\ZZZ H$-modules
(here $\D$ stands for the triangulation of $X$).
There is a chain homotopy equivalence
\bq
\chi_\D:C_*^\D(\wi X)\to C_*^s(\wi X)
\end{equation}
of chain complexes over $\ZZZ H$, functorial  up to
chain homotopy (see \cite{spanier}, Ch. 4, \S 4.)
The next lemma is a corollary  of
the well known result on the combinatorial invariance of the
Whitehead torsion (see \cite{milnWT}).
\bele\label{l:chap}
Let $X$ be a $\smo$ closed manifold.
Let $\D_1, \D_2$
be $C^1$-triangulations of $X$.
Then
\bq
\chi_{\D_1}^{-1}\circ\chi_{\D_2}  :
C_*^{\D_2} (\wi X)
\to
C_*^{\D_1} (\wi X)
\end{equation}
is a simple \heq~
of finitely generated free based chain
complexes over $\ZZZ H$ (that is, its torsion 
$\tau(
\chi_{\D_1}^{-1}\circ\chi_{\D_2})$
vanishes in $\Wh(H)=K_1(\ZZZ H)/(\pm H)$).
 $\qs$
\enle

Return now to the Morse-Novikov theory.
Let $(M,f,v)$ be a $MN$-triple.
In the previous discussion let us set $X=M$
and consider the infinite cyclic covering 
$\bar M \to M$
corresponding to $\xi=\xi(f)$.
Let $\D$ be a triangulation of $M$
and define a chain map

\bq
\Phi_\D=\chi^{-1}_\D \circ \Phi(M,f,v)
:
C_*(f,v)\to
 C_*^\D(\bar M)\tens{L}\wh L
\end{equation}
The chain complexes
$C_*(f,v)$ and 
$C_*^\D(\bar M)\tens{L}\wh L$
are 
finitely generated 
 free complexes over $\wh L$.
Choosing a lift of a 
subset $S(f)\sbs M$ to $\bar M$ and  lifts
of all the simplices of $\D$
to $\bar M$
we obtain free bases in both complexes.
With respect to any such choice of bases
we obtain the torsion 
$\tau(\Phi_\D)
\in \ove K_1(\wh L)$.
The ambiguity in the choice of bases leads to
 multiplication by $\pm t^n$,
so that the image of $\tau(\Phi_\D)$
in $\ove K_1(\wh L)/T$
does not depend on these choices and, 
moreover, does not depend on the choice of $\D$ by Lemma 
\ref{l:chap}.

\bede\label{d:w}
The 
element 
$\det(\tau(\Phi_\D))\in W$
will be denoted by $w(M,f,v)$.
\end{defi}
Our  second 
main result -- Theorem B -- says
that $w(M,f,v)$ is equal to the inverse Lefschetz zeta function 
(Subsection \ref{su:eta_zeta}) of the flow
generated by
$(-v)$. For convenience of notation we shall abbreviate 
$(\zeta_L(-v))^{-1}$
as $\zeta(v)$.

\begin{B}\label{t:B}
For every $v\in \GG(f)$
we have
\bq
w(M,f,v)=  \zeta(v).
\end{equation}
\end{B}

Now we shall say some words about the proofs. Theorem A
is proved by analysing step by step the construction from
\cite{patou}.
The arguments which we use here to prove the functoriality
are based on noetherian commutative algebra.
The proof of Theorem B is based on our work
\cite{pastpet}.
In this paper we have introduced a class $\GG\TT_0(f)$
of $f$-gradients 
(where $f:M\to S^1$
is a given Morse function)
which is $C^0$-open-and-dense in the set 
$\GG\TT(f)$
of all $f$-gradients satisfying \TA.
For every
$v\in 
\GG\TT_0(f)$
the dynamics of the gradient flow generated  by $v$
has many remarkable properties; in particular
the Novikov incidence coefficients
are rational functions.
A slight modification of the definition
of the class
$\GG\TT_0(f)$
(\cite{pafest})
leads to the subset
$\GG_0(f)\sbs
\GG(f)$
which is $C^0$-open-and-dense
in $\GG(f)$;
for every 
$v
\in
\GG\TT_0(f)$
the Novikov incidence coefficients
are rational functions
and the formula \rref{f:tors_p}
holds.
The case of general gradient $v$
as required by Theorem B
is done by a $C^0$-perturbation argument
which makes the reduction to the $C^0$-generic case.
Actually we do not know if both sides of
\rref{f:tors_p}
are invariant under such perturbation.
But note that both these expressions are power series in $t$.
We show that each {\it finite part}
of such
series is invariant \wrt~such perturbation, and this finishes the proof.

\subsection{Related work}
\label{su:rel_work}

In the paper  \cite{farran}
M.Farber and A.Ranicki gave an alternative construction of a chain
complex generated by critical points of a circle-valued
Morse map $f:M\to S^1$, and proved that this complex computes
the completed homology of the infinite cyclic covering.
(The paper \cite{farran} contains actually a more general 
construction of a chain complex over Novikov completions of fundamental
group.) In the preprint \cite{ranXCIX}
A.Ranicki identified this chain complex in the particular case of
the infinite cyclic covering 
with the chain complex
from \cite{pastpet}.
The results of the paper \cite{pafest} were generalized
in \cite{pawitt} to the non-abelian case. 
The work of D.Sch\"utz  \cite{schuetzC}
contains a generalization of the results of
\cite{pawitt}
to the case of irrational forms.
One of the ingredients in the work of D.Sch\"utz,
which is quite new to the Morse-Novikov theory
is the application of Hochshild homology techniques,
which is related to the work  \cite{gnXCIV} by R.Geoghegan and A.Nicas.
The results of \cite{schuetzC}
pertained to the $C^0$-generic case;
in a recent preprint \cite{schuetzCI}
D.Sch\"utz considers the general case. He constructs, in
particular, a new chain equivalence between 
Novikov complex and the completed simplicial chain complex.

\subsection{Acknowledgements}
\lb{su:ackno}

It is a pleasure for me to express here my gratitude 
to V.Turaev  first of all for 
his suggestion in 1998 of the precise
form of the functoriality in Theorem A,
but not the less 
for many useful  discussions and critics, 
which influenced strongly the initial plan of the paper.

The final step of the work  was done  during my stay 
at ETHZ in the fall of 2000. I am grateful to ETHZ for 
the excellent working conditions.
Many thanks to E.Zehnder and D.Salamon for the warm hospitality
during my stay.

\section{Preliminaries on chain complexes}
\label{s:cha_com}

In this section we work in the category of chain complexes
of left  $R$-modules, where $R$ is a ring.
We consider in this paper only chain complexes,
concentrated in positive degrees.
We shall often omit the adjective "chain", so that 
{\it complex }~ means "chain complex", 
{\it map }~ means "chain map", 
{\it homotopy } means "chain homotopy" etc. 
A chain map $f:C_*\to D_*$
is called {\it homology equivalence}
if it induces an isomorphism in homology.
We shall wrtie $f\sim g$ if $f$ is homotopic to $g$.

\subsection{Models}
\label{su:models}

Let $R$ be a ring (not necessarily commutative).
\bede
Let $X_*, Y_*$
be complexes over $R$.
A map
$f:X_*\to Y_*$
is called {\it model for $Y_*$}
if $f$ is a \heeq~
and $X_*$
is free.
\end{defi}
The next
proposition states that for a given $Y_*$
its model is essentially unique.
\bepr\label{p:models}
Let
$f:X_*\to Y_*,~
f':X'_*\to Y_*$
be models for $Y_*$.
Then there is a \heq~$\mu:X_*'\to X_*$,
\sut
~ 
$f\circ\mu \sim f'$. 
The homotopy class of $\mu$
is uniquely determined by the homotopy classes of $f$ and $f'$.
\enpr
This proposition follows from  the next one.

\bepr\label{p:mod_gen}
Let $\a:X_*\to Y_*, ~  \b:Z_*\to Y_*$
be chain maps, and assume that $X_*$
is free and $\b$
is a \heeq.
Then there is a map
$\g:X_*\to Z_*$, \sut~ $\b\circ\g\sim\a$.
The homotopy class of $\g$
is uniquely determined by the homotopy classes of $\a$ and $\b$.
\enpr

\Prf
We shall prove the existence of $\g$, the homotopy uniqueness is 
proved similarly.
Adding to $Z_*$ a contractible chain complex if necessary,
we can assume  that the map $\b$ is epimorphic.
Using induction on degree we shall construct a map
$\g$ satisfying
$\b\g(x)=\a(x)$.
Assume that $\g$
is already defined on every $X_i$
with $i<k$
and satisfies the equation $\b\g(x)=\a(x)$
for $\deg x<k$.
It suffices to construct for every free generator
$e_k$ of $X_k$
an element
$x\in Z_k$, \sut
$\b(x)=\a(e_k)$ and $\quad \pr x=\g(\pr e_k)$.
Choose any $x_0\in Z_k$, \sut~
$\b(x_0)=\a(e_k)$.
Consider the element
$y=\pr x_0-\g(\pr e_k)$.
It is a cycle of the complex
$Z'_*=\Ker\b$.
Since $\b$ is a \heeq,
the complex $Z'_*$
is acyclic, therefore, there is
$z\in Z'_k$, \sut
~$\pr z = y$.
Now set 
$x=x_0-z$
and the proof is over. $\qs$

\subsection{Filtrations and adjoint complexes}
\label{su:fil_adj}
We recall here briefly the material of \cite{patou}, \S 3.A.
We shall call {\it filtration}
of a chain complex
$C_*$
a sequence of subcomplexes $C_*^{(i)}, -1\leq i$,
\sut~
\bq
0= C_*^{(-1)}
\sbs
C_*^{(0)}
\sbs
...
\quad \mbox{ and }
\quad
\cup_i
C_*^{(i)}
=
C_*
\end{equation}
A filtration $C_*^{(i)}$
is {\it good}
if
$H_k(
C_*^{(i)}
/
C_*^{(i-1)})
=
0
\mbox{ for } k\not= i$.
For a filtration
$\{
C_*^{(i)}
\}$
of a complex $C_*$
set
$C_n^{gr}=
H_n(
C_*^{(n)}
/
C_*^{(n-1)}
)$,
and let
$\pr_n:C_n^{gr}\to C_{n-1}^{gr}$
be the boundary operator of the exact sequence of the triple
$(
C_*^{(n)},
C_*^{(n-1)},
C_*^{(n-2)})
$.
Then $C_*^{gr}$
endowed with the boundary operator $\pr_n$
is a complex, which will be called
{\it adjoint}
to $C_*$.

{\bf Example. }
Let $D_*$
be any complex.
The filtration
\bq
D_*^{(i)}
=
\{
0
\leftarrow
 D_0 {  \leftarrow
}  {... } \leftarrow
D_i   \leftarrow
0 { \leftarrow
}...
\}
\end{equation}
is called {\it trivial}.
This is obviously a good filtration and
$D_*^{gr}=D_*$.
The proof of the next lemma is in standard diagram chasing.
\bele[\cite{patou}, Lemma 3.2]\label{l:cellular}
Let
$C_*, D_*$
be complexes.
Assume that $C_*$
is endowed with a good filtration, and that
$D_*$ is a free complex, endowed
with the trivial filtration.
Let
$\phi:D_*\to C_*^{gr}$
be a map.
Then there exists a map
$f:D_*\to C_*$
preserving filtrations and inducing
the map $\phi$
in the adjoint complexes.
The map $f$
is unique up to a homotopy preserving filtrations. $\qs$
\enle

\bede
A good filtration
$\{
C_*^{(i)}
\}
$
of a complex $C_*$
is called {\it nice}
if every module
$H_n(
C_*^{(n)}
,
C_*^{(n-1)}
)$
is a free $R$-module.     \end{defi}

\beco\label{c:heeq}
For a nice filtration
$\{C_*^{(i)}\}
$
of a complex
$C_*$
there exists a \heeq~
$C_*^{gr} \to C_*$
functorial up to homotopy
in the category of nicely filtered complexes. 
If $C_*$ is a complex of free $R$-modules, 
this  \heeq~is a \heq.
$\qs$
\enco

\subsection{Strings and inverse limits}
\label{su:strings}

An infinite sequence
\bq
\CC
=
\{C_*^0
\leftarrow
   C_*^1
\leftarrow
{...}
\}
\end{equation}
of chain epimorphisms
is called {\it string}.

A {\it map of strings}
$\CC\rTo^h \DD$
is a diagram of the following type:
\begin{diagram}[size=2em,LaTeXeqno]
C_*^0 & \lTo & C_*^1 & \lTo &..... &\lTo
& C_*^n & \lTo^{p_n} & C_*^{n+1}
 ...\\
\dTo_{h_0} & &  \dTo_{h_1}& &...&
 & \dTo_{h_n}& &\dTo_{h_{n+1}}\lb{f:mapstrings} \\
D_*^0  & \lTo & D_*^1  & \lTo & ..... &\lTo
  &  D_*^n & \lTo^{q_n} & D_*^{n+1}...
   \\
\end{diagram}
where $h_i$ are chain maps and 
all the squares are homotopy commutative.

A {\it strict map of strings}
is a map of strings where all the squares
in (\ref{f:mapstrings})
are commutative.
For a string $\CC$
the chain complex
$\varprojlim C_*^i$
will be  denoted by
$|\CC|_*$ and called {\it inverse limit of $\CC$}.

A strict map of strings induces obviously a chain map
of their inverse limits. 
The aim of the next proposition is to generalize this 
property to the case of arbitrary maps of strings.
This proposition (and its proof) is  quite  close to  Proposition  3.7 of \cite{patou},  
so we shall give only a sketch of the proof.

\bepr\label{p:map_strings}
Let
$h:\AA\to \BB$
be a map of strings, $h=\{ h_i\}$.
Assume that $|\AA|_*$
is a chain complex of free $R$-modules.
Then there is a chain map
$\HH:|\AA|_*\to |\BB|_*$, \sut
~for every $k$
the following diagram is homotopy commutative:

\begin{diagram}[size=2em,LaTeXeqno]
|\AA|_* & \rTo^{\HH} &   |\BB|_*   \\
\dTo &   &\dTo  \qquad\lb{f:dopred_usl}\\
A_*^k   & \rTo^{h_k} & B_*^k \\
\end{diagram}

(where the vertical arrows are natural projections)
\enpr

{\it Sketch of proof. \quad}
For $k\geq 0$
let $Z_*^k$
be the chain cylinder of the map
$h_k:A_*^k\to B_*^k$.
Using homotopy commutativity of
 the squares in (\ref{f:mapstrings})
one defines a string $\ZZ
=
\{
Z_*^0
\leftarrow
Z_*^1
\leftarrow
...
\}
$
together with 
two strict maps of strings
\bq
\AA\rTo^\a \ZZ, \quad \BB\rTo^\b \ZZ
\end{equation}
\sut~ 
~the corresponding maps
$\a_k:A_*^k\to Z_*^k, \b_k:B_*^k\to Z_*^  k$
are the standard inclusions of the source, resp.
the target to the cylinder.
Form the corresponding maps of the inverse limits:
\bq
|\AA|_*\rTo^{\wi\a} |\ZZ|_*, \quad |\BB|_*\rTo^{\wi\b} |\ZZ|_*
\end{equation}
The map
$\wi\b$ is a \heeq.
(Indeed, 
by the very  definition every
string satisfies the Mittag-Leffler condition
(\cite{massey}, definition A16). Therefore for every $n$ we have 
$\lim^1 B_n^i=0,~ \lim^1 Z_n^i=0$. In view of 
\cite{massey}, Th. A.19 
it suffices to recall  that each 
$\b_k$ is a homology equivalence.)
Apply now Proposition \ref{p:mod_gen}
to obtain a chain map
$\HH:|\AA|_*\to |\BB|_*$,
 \sut~
$\wi\b\circ\HH\sim\wi\a$.
The commutativity of the diagram
(\ref{f:dopred_usl})
is now easy to check
using Proposition \ref{p:mod_gen}. $\qs$

It is natural to ask whether the chain homotopy class of
$\HH$
is uniquely determined by the condition
(\ref{f:dopred_usl}).
The  next two subsections
provide an answer to this question for the particular case
of the commutative noetherian base ring.

\subsection{A bit of noetherian homological algebra}
\label{su:noet}

We interrupt for a moment our study of
 strings and their inverse limits
in order to prove a proposition about chain complexes
over commutative power series rings.
In this section $A$ is a commutative noetherian ring,
$R=A[[t]]$, and $R_n=A[[t]]/t^n$ (where $n\geq 0$).
A finite chain complex 
of free $R$-modules
will be called {\it homotopy finitely generated }
\label{d:hom_fin}
if it is chain equivalent to
a finite chain complex of finitely generated 
free $R$-modules.

\bepr\label{p:noet}
Let $C_*, D_*$
be homotopy finitely generated  complexes
over $R$. 
Let 
$f,g:C_*\to D_*$
be chain maps.
Let
$C_*^{(n)}=C_*/t^n C_*, \quad D_*^{(n)}=D_*/t^n D_*$,
and
 $f_n=f/t^n, g_n=g/t^n: C_*^{(n)} \to D_*^{(n)}$
be the corresponding quotient maps.
Assume that $f_n\sim g_n$ for every $n\in\NNN$.
Then $f\sim g$.
\enpr
\Prf
It suffices to prove the proposition for  the particular case
when both $C_*, D_*$
are finitely generated. In this assumption, 
let $H$ be the set of
{\it homotopy classes }
of maps
$C_*\to D_*$;
then $H$
is a finitely generated $R$-module.
Note that the maps $f_n, g_n$
are homotopic if and only if the map
$f-g:C_*\to D_*$
is homotopic to a map divisible by $t^n$.
The proof is finished by the following
lemma.
\bele
Let $H$ be a finitely generated
$R$-module. Let $x\in H$. Assume that $x$ is divisible by $t^n$
for every $n$. Then $x=0$. 
\enle
\Prf
Let us first consider the case, when $H$
has no $t$-torsion (that is, $tx=0$ 
implies $x=0$ for $x\in H$). Let $N\sbs H$
be the submodule of elements, divisible by $t^n$
for any $n$. Then $N$  satisfies
$tN=N$. 
Indeed, let $x\in N$.
Then for every $n\in\NNN$
we have
$x=t^n y_n$.
Moreover, there is only one $y_n$
satisfying this relation, and $y_n=ty_{n+1}$
(this follows since $N$ has no $t$-torsion).
Now $y_1$ is obviously divisible by all powers of $t$.
Since $t$ is in the Jacobson radical of $R$, the Nakayama's lemma
implies $N=0$.

Let us now consider the case of arbitrary $H$.
Let $x$ be an element divisible by every  power of $t$.
Let $T\sbs H$
be the submodule of all $t$-torsion elements, i.e.
$T=\{x~|~\exists n\in\NNN, t^nx=0\}$.
Set $H'=H/T$. Applying to $H'$
the reasoning above, we conclude that $x\in T$.
Choose some $k\in\NNN$, such that $t^k T=0$.
Since $x=t^ku$ (with $u$ necessarily in $T$), we have $x=0$.
$\qs$

\subsection{Strings and inverse limits: part 2}
\label{su:strings2}
As in the previous section let $A$ be a noetherian commutative ring
and $R=A[[t]]$.
Let $C_*$ be
a free finitely generated
complex  of $R$-modules.
The string 
\bq
\CC= \{C_*/tC_*\lTo C_*/t^2C_* {\lTo }... \lTo C_*/t^n C_* {\lTo} ... \}
\end{equation}
 will be called 
{\it the special string corresponding to $C_*$}.
Note that since $C_*$ is free and finitely generated,
there is a natural isomorphism
$|\CC|_*\approx C_*$.

\bepr
In the hypotheses of Proposition
\ref{p:map_strings}
assume moreover that $\AA$ and $\BB$
are special strings. 
There is only one  (up to homotopy)
map 
$\HH:|\AA|_*\to |\BB|_*$
\sut~all the squares 
(\ref{f:dopred_usl}) 
are homotopy commutative.
\enpr
\Prf
Let 
$A_*, B_*$
be the free finitely generated complexes over $P$,
generating $\AA$, resp. $\BB$.
If 
$  \HH_1, \HH_2  :|\AA|_*\to |\BB|_*$
both make all the squares 
(\ref{f:dopred_usl})
commutative, then
the maps 
\bq
\HH_1/t^k, \HH_2/t^k: |\AA|_*/t^k|\AA|_*\to |\BB|_*/t^k|\BB|_*
\end{equation}
are homotopic for every $k$.
Apply now Proposition \ref{p:noet}
and the proof is over. $\qs$

Therefore any
map
$h:\AA\to\BB$
of special strings determines a chain map
$|\AA|_*\to|\BB|_*$
which is well defined up to homotopy;
this map will be denoted $|h|$.
The next corollary is immediate.

\beco\lb{co:square}
Let
$\AA=(A^i_*),~ \BB =(B_*^i),~ \CC=(C_*^i), ~ \DD=(D_*^i) $
be special strings, and

\begin{diagram}[LaTeXeqno]
\AA & \rTo^{\a} & \BB\\
\dTo_\psi &   &\dTo_{\phi}   \lb{f:ssquare}\\
\CC     & \rTo^{\b} & \DD\\
\end{diagram}

be a square of strings, which is chain homotopy commutative on each 
finite level, i.e.
for every $n$ the chain maps
$\phi_n\circ\a_n,~\b_n\circ\psi_n:A_*^n\to D_*^n$
are homotopic. 
Then the following square is chain homotopy commutative.

\begin{diagram}[LaTeXeqno]
|\AA_*| & \rTo^{|\a|} & |\BB_*|\\
\dTo_{|\psi|} &   &\dTo_{|\phi|}   \lb{f:sssquare}\\
|\CC_*|     & \rTo^{|\b|} & |\DD_*|\\
\end{diagram}
$\qs$
\enco

\section{Proof of theorem A}
\label{s:proof_a}

We start by giving more details on  the Morse complex and the 
Morse-Novikov complex
(Subsections \ref{su:more_morse}, \ref{su:nov_basics}).
In Subsection \ref{su:morse_wn}
we recall some points from \cite{patou}
which will be useful in the sequel. 
The construction of a functorial chain equivalence
$\Phi$ is done in Subsection \ref{su:equi_fi} after a technical 
Subsection \ref{su:spec_t_ord}.
The end of the proof of Theorem A is in Subsection \ref{su:functo}.

Now some words about terminology.
The term
{\it $f$-gradient}
means {\it a gradient-like vector field for $f$}
(see \cite{milnhcob}, \S 3 for  definition).
For a $C^1$ vector field $v$ on a \ma~$M$
the symbol
$\g(x,t;v)$
 denote  
the value at $t$ of the integral curve of $v$
which passes through $x$ at $t=0$.

\subsection{About Morse complexes}
\label{su:more_morse}

Let 
$\fcob$
be a Morse function on a cobordism $W$,
$v$ be an $f$-gradient.
The set of critical points of $f$
will be denoted by $S(f)$,
the set of critical points
of $f$ of index $k$
is denoted by
$S_k(f)$. For  $p\in S(f)$
let $D(p,v)$ denote the stable manifold of $p$ \wrt~$v$:
\bqq
D(p,v)=\{x\in W~|~\g(x,t;v)\rTo_{t\to\infty} p\}
\end{equation*}
We assume that $v$ satisfies {\it transeversality condition}
that is 
for every $p,q\in S(f)$
we have
\bqq
D(p,v)\pitchfork
D(q, -v).
\end{equation*}
For every $p\in S(f)$ 
choose an orientation
of $D(p,v)$.
To this data one associates a chain complex
$C_*$
of free abelian groups as follows.
By definition $C_k$
is a free abelian group
generated by $S_k(f)$.
To define the boundary operator,
let $p\in S_k( f), q\in S_{k-1}( f)$,
and let
$\G(p,q;v)$
be the set of all orbits of $v$, joining $p$ to $q$.
The choice of orientations allows to attribute to each orbit
$\g\in
\G(p,q;v)$
a sign $\ve(\g)\in \{-1, 1\}$.
Set $n(p,q;v)
=
\sum_{\g\in \G(p,q;v)} \ve(\g)$
and define a \ho~
$\pr_k:C_k\to C_{k-1}$
by
\bq\label{f:bou_ma}
\pr_k p=\sum_{q\in C_{k-1}( f)} n(p,q;v)\cdot q
\end{equation}
One can check that
$\pr_{k-1}\circ\pr_k=0$.
The resulting chain complex
is called the {\it Morse complex}
or the {\it Morse-Thom-Smale-Witten complex}.
We shall denote it
$C_*^\MM(f,v)$,
or simply
$C_*(f,v)$
if \noconf.
We shall now outline  the construction of a chain
\heq
\bq\label{f:nat_heq}
C_*(f,v)\to C_*^s( W, {\pr_0 W})
\end{equation}
(following \cite{patou}, Appendix,
where the reader 
will find the  details).
A Morse function
$\phi:W\to [a,b]$
is called {\it ordered}
if $\phi(x)<\phi(y)$
whenever $x,y\in S(\phi)$
and
$\ind x<\ind y$.
By a standard application
 of the rearrangement procedure, see \cite{milnhcob}, \S 4)
there is an ordered Morse function
$\phi:W\to [a,b]$
\sut
~$v$ is also a $\phi$-gradient.
Let
$a=a_0<a_1< ... <a_m<a_{m+1}=b$
be an {\it ordering sequence } for $\phi$,
that is, every
$a_i$
is a regular value for $\phi$
and
$S_i(\phi)\sbs \phi^{-1}([a_i,a_{i+1}])$
for every $i$
(here $m$ stands for the dimension of $M$).
Let
$
W^{(i)}
=\phi^{-1}([a_0,a_{i+1}])$
and consider the filtration
of the pair
$( W, \pr_0 W)$
by the pairs
$( W^{(i)}
 , \pr_0 W)$.
The standard Morse theory (see \cite{milnmt})
implies that (up to homotopy equivalence) $W^{(i)}$
is obtained from
$W^{(i-1)}$
by attaching cells of dimension $i$, therefore
the corresponding filtration
in the singular chain complex
$C_*^s( W, {\pr_0 W})$
is {\it nice}.
The corresponding adjoint complex $D_*$
(see \ref{su:fil_adj})
is 
freely generated (as  abelian group) in dimension $k$ by
$S_k(f)$.
One can show that the boundary operator in $D_*$
is given by the formula
(\ref{f:bou_ma})
(the argument is the same as in \cite{milnhcob}, Corollary 7.3).
Thus the Morse complex
$C_*^\MM(f,v)$
is identified
with 
$D_*$.
The existence of the natural \heq~
(\ref{f:nat_heq})
follows now from 
\ref{c:heeq}.

\subsection{Basics of Morse-Novikov theory}
\label{su:nov_basics}

Let $(M,f,v)$ be a Morse-Novikov triple, as defined in \ref{su:stat_res}.
As in the case of real-valued Morse functions 
we denote by  $D(p,v)$ the stable manifold of $p$ \wrt~$v$.
Choose an orientation of $D(p,v)$ for every $p\in S(f)$.
To this data one associates a chain complex
$C_*$
of free finitely generated
$\wh L$-modules as follows.
Let $ F:\bar M\to \RRR$
be a lift of $f:M\to S^1$.
Let $S_k(F)$denote the set of critical points of $F$ of index $k$.
Consider the set $C_k$
of all formal linear combinations $\l$ of 
elements of $S_k( F)$, \sut~
for any
$c\in\RRR$
there is only a finite set of points in $\supp\l$
which lie above $c$.
It is easy to see that
$C_*$
is a free $\wh L$-module
(any lift of the set $S_k(f)$
to $\bar M$
provides a family of free $\wh L$-generators for $C_k$).
To define the boundary operators, let
$p\in S_k(F), q\in S_{k-1}(F)$
and denote by
$\G(p,q;v)$
the set
of all orbits of $v$ in $\bar M$,  joining $p$ to $q$.
The transversality condition
implies that
$\G(p,q;v)$
is finite.
The choice of orientations allows to attribute
to each flow line
$\g\in
\G(p,q;v)$
a sign $\ve(\g)\in \{-1, 1\}$.
Set
$n(p,q;v)
=
\Sigma_{\g\in
\G(p,q;v)}
\ve(\g)$
and define a \ho
~$\pr_k:C_k\to C_{k-1}$
by
\bq\label{f:bou_nov}
\pr_k p=\sum_{q\in S_{k-1}(f)} n(p,q;v) \cdot q 
\end{equation}

One can check that $\pr_{k-1}\circ \pr_k=0$.
The resulting chain complex is called the
{\it Novikov complex}.
We denote it
$C_*^\NN(M,f,v)$
or simply $C_*(f,v)$
if \noconf.

\subsection{Morse complexes of finite pieces of the cyclic covering}
\label{su:morse_wn}

Working with the terminology of \ref{su:nov_basics},
choose a regular value $\l$
of $F:\bar M\to\RRR$ and set
$V=F^{-1}(\l)$.
For $\a\in\RRR$
set
$V_\a=F^{-1}(\a)$.
Set
$$
W=F^{-1}\([\l-1, \l]\), \quad
V^-=F^{-1}\(]-\infty,\l]\).
$$
Thus the  cobordism $W$ is  the result of
cutting of $M$ along $V$.
The structure group
of the covering $\bar M\to M$
is isomorphic to $\ZZZ$ and we choose the generator $t$
of this group so that $tV_\a =V_{\a-1}$.

We have
\bqq
\bar M=
\bigcup_{s\in\ZZZ} t^s W,    
\quad  \mx{     ~with ~ }\quad     t^{s+1} W \cap   t^s W  = V_{\l-s-1}.
\end{equation*}
For any $k\in\ZZZ$ the map $t^k:V_{\l+k}\to V$
is a diffeomorphism.
For $n\in\ZZZ, n\geq 1$
let
\bq
W_{n}=\bigcup_{0\leq s\leq n-1} t^{s}W=
F^{-1}([\l-n, \l])
\end{equation}
so that in particular
$W=W_1$.
The Morse complex
\bq
\MM_*(n) = C_*^\MM(F|W_n,v|W_n)
\end{equation}
of the function $F|W_n:W_n\to[\l-n,\l]$ 
is a chain
complex of free abelian groups.
But it has some more structure,
coming from the $\ZZZ$-action  
on $\bar M$. Set
$$
\wh L_-=\ZZZ[[t]], \quad L_n=\ZZZ[t]/t^n, \quad  L_-=\ZZZ[t]
 $$
It  is easy to show, using the $\ZZZ$-invariance of $v$, that
$\MM_*(n)$
is a free chain complex over $L_n$.   
In particular, $\MM_*(n)$
is a chain complex of $\wh L_-$-modules.
It is clear that there is a natural isomorphism
\bq\label{f:n_n+1}
\MM_*(n+1)\tens{\wh L} L_n \approx \MM_*(n)
\end{equation}
Any lift to $W$ of $S(f)$
provides a free $L_n$-base for $\MM_*(n)$ and 
the isomorphsism
(\ref{f:n_n+1})
preserves these bases.
Introduce now the string
\bq\label{f:string_m}
\MM
=
\{
\MM_*(1)
\lTo^{\pi_2}
\MM_*(2)
{\lTo^{\pi_3}}
{...}
\lTo^{\pi_n}
\MM_*(n)
{\lTo}
{...}
\}
\end{equation}
where $\pi_n$
is the projection
$
\MM_*(n)
\to
\MM_*(n)
\tens{ L}
 L_{n-1}
$.
The inverse limit
$|\MM|_*$
is a free
$\wh L_-$-complex,
which will be denoted by
$C_*^-(f,v;\l)$.
It is clear that there is a base preserving isomorphism
\bq
C_*^\NN(M,f,v)
\approx
C_*^-(f,v;\l)\tens{\wh L_-}\wh L.
\end{equation}
Therefore the  Novikov complex can be
reconstructed from the Morse complexes
$\MM_*(n)$
if we take into account the structure of
$ L_n$-modules on these complexes.

The natural chain equivalence
$\MM_*(n)
\to
C_*^s( W_n, {\pr_0 W_n})
$
can be refined so as to respect the
$ L_n$-structure.
Now we shall recall this construction
(from \cite{patou}, \S 5)
in more details.
The chain equivalence
$\MM_*(n)
\to
C_*^s( W_n, {\pr_0 W_n})
$
is constructed using an ordered
Morse function on $W_n$.
It turns out that if we choose an ordered Morse function
on $W_n$
which is well fitted to the action of $\ZZZ$
on $\bar M$,
the chain \heq
~will respect the $ L_n$-structure.
\bede
An ordered Morse function
$\phi:W_n\to[\a,\b]$
is called
$t$-ordered
if for every
$x\in W_n$
we have
\bq\label{f:t_ord}
tx\in W_n\Rightarrow
\phi(tx)<\phi(x).
\end{equation}
\end{defi}

\bepr[\cite{patou}, Lemma 5.1]\lb{p:t_ord}
There is a $t$-ordered Morse function $\phi$
on $W_n$, \sut~
$v$ is a $\phi$-gradient. $\qs$
\enpr

Let us see how the proposition applies.
Let $\phi$
be any $t$-ordered Morse function
on $W_n$
\sut~$v$ is a $\phi$-gradient.
Then $\phi$ induces a filtration 
$W^{(i)}$
on $W$ as explained in \ref{su:more_morse}.
Set 
\bq\label{f:x_i}
X^{(i)}=W^{(i)}\cup t^n V^-.
\end{equation}
The chain complexes
$C_*^s(X^{(i)}, t^n V^-),
 ~0\leq i \leq \dim M$
form a nice filtration of
$C_*^s(V^-, t^n V^-)$
(with $L_n$ as the base ring).
The adjoint complex is
$\MM_*(n)$
and thus we obtain a chain \heq
\bq\lb{f:homeq_J}
J_n:\MM_*(n)\to
C_*^s( V^-, {t^n V^-})
\end{equation}
of chain complexes over $\wh L_n$.
One can show that the chain homotopy class
of $J_n$
does not depend on the particular choice of $\phi$
(\cite{patou}, p. 324).

\subsection{Special $t$-ordered functions}
\label{su:spec_t_ord}

In view of  later applications it will be useful to choose 
the $t$-ordered function 
$\phi$ 
so that
the terms of the corresponding filtration
of $W_n$
are very close to the descending discs of $v$. 
To make this  precise,
we introduce first a  definition.
Set 
\bq\label{f:indl_s}
D(\indl i;v)
=
\bigcup_{\ind p\leq i} D(p,v)
\end{equation}
where
$D(p,v)$
is the stable \ma~in $W_n$
of the point
$p\in S(F)$,
and let
$U_i$
be any \nei~
in $W_n$
of
$D(\indl i;v)\cup\pr_0 W_n$.
Let $m=\dim M$.

\bepr\label{p:thin_fi}
There is a $t$-ordered Morse function
$g:W_n\to \RRR$,
\sut
~$v$ is a $g$-gradient,
and an ordering sequence
$a_0<... <a_{m+1}$
for $g$, \sut
~the corresponding filtration
$\{W^{(i)}\}$
satisfies $W^{(i)}
\sbs U_i$
for every $i$.
\enpr
\Prf
Start with any
$t$-ordered Morse function
$h:W_n\to[\a,\b]$.
We can assume that
for an $\e>0$ the function
$dh(x)(v(x))$
is constant in
$
h^{-1}([\b-\e,\b])
\cup
h^{-1}([\a,\a+\e])
$.
Choose $\e$ so small that
$h^{-1}([\a,\a+\e])
\sbs U_i
$ for every $i$.
Let
$\mu:[\a,\b]\to[0,1]
$
be a $\smo$
function, \sut
\bq
\left\{
\begin{aligned}
\mu(y)=0& \mbox{  \quad  for\quad }
y\in [\a,\a+\e/2]
\mbox{\quad  or \quad  }
y\in [\b-\e/2,\b]
\\
\mu(y)=1& \mbox{\quad  for\quad  }
y\in [\a+\e,\b-\e]
\end{aligned} \right.
\end{equation}

Set $w(x)=\mu(h(x))v(x)$
and let $\Phi_\tau$
denote the one-parametric group of diffeomorphisms
of $W_n$, induced by $w$.
I claim that for $T>0$
sufficiently large the function
$g=h\circ\Phi_T$
satisfies the conclusion of our Proposition.
It is quite obvious that $v$ is an $g$-gradient.
Next we shall check that $g$ is $t$-ordered.
Let
$
\g_1(\tau)=\g(x,\tau, w),
\g_2(\tau)=\g(tx,\tau, w)$.
We have to prove
\bq\label{f:nervo_t_ord}
h(\g_1(T))>
h(\g_2(T))
\end{equation}

We shall show that (\ref{f:nervo_t_ord})
holds for $x$ satisfying the following condition:
both $x$ and $tx$ are not in $\supp(v-w)$, and both $\g_1, \g_2$
reach $\pr_1 W$
at some moment.
(The case of arbitrary $x$ is similar.)
Let $t_1$ (resp. $t_2$)
be the moment when $\g_1$, (resp. $\g_2$)
reaches $h^{-1}(\b-\e)$.
Note that
$\g_2(\tau) = t\g_1(\tau)$
for
$\tau\in [0,t_1]$.
Therefore $t_2>t_1$
and (\ref{f:nervo_t_ord})
holds if $T\leq t_1$.
It holds also for $T\in [t_1, t_2]$, since for these values of $T$
we have
$h(\g_2(T))
\leq
\b-\e
\leq
h(\g_1(T))$.
Finally, for $T\geq t_2$
both
$\g_1(T), \g_2(T)$
are in
$h^{-1}([\b-\e,\b])$
and
$h(\g_1(T))=h(\g_2(T+t_2-t_1))>h(\g_2(T))$.

Any
ordering sequence
$\a_0<\a_1<... <\a_{m+1}$
 for $h$
is also an ordering sequence for $g$.
A
standard "gradient descent" argument
shows that
$g^{-1}([\a_0, \a_{i+1}])\sbs U_i$
if only $T$ is sufficiently large
and $\a_{m+1}<\b-\e$.
$\qs$

\subsection{Construction of the chain \heq~$\Phi$}
\label{su:equi_fi}

Introduce the following notation:
\bq
S_*= 
C_*^s( V^-)\tens{ L_-}\wh L_-
,\quad
S_*(n)=
C_*^s( V^-, t^n V^-)
\end{equation}
Then $S_*$
is a free 
complex over $\wh L_-$,
and $S_*(n)$
is free over $ L_n$.
The chain complex $S_*$
is infinitely generated, but it is 
homotopy finitely generated (in the sense of the definition on the page
\pageref{d:hom_fin}).
Indeed, choose a $C^1$-triangulation $\D$
of $M$, \sut~
$V$ is a simplicial subcomplex. Then $ V^-$
obtains a  triangulation, invariant \wrt~
the action of $t$, 
so that the simplicial chain complex 
$
C_*^\D( V^-)
$
is a free finitely generated complex over $ L_-$.
Let
$D_*
=
C_*^\D( V^-)
\tens{ L_-}\wh L_-$.
the natural chain equivalence
$
C_*^\D( V^-)
\to
C_*^s(V^-)$
being tensored with $\wh L_-$
gives then the homotopy equivalence sought.

\bepr\label{p:mu}
There is a unique (up to homotopy) chain \heq~
$
C_*^-(f,v;\l)
\rTo^{\mu}
S_*
$, \sut~
for every $n\geq 1$ 
~the following diagram is homotopy commutative:
\begin{diagram}[size=1.8em,LaTeXeqno]
 C_*^-(f,v;\l) & \rTo^{~\mu} & S_* \\
\dTo &            &\dTo   \quad\lb{f:cond_n}\\
\MM_*(n)   & \rTo^{J_{n}} & S_*(n)\\
\end{diagram}
(where the vertical arrows are natural projections, and $J_n$
is the chain homotopy equivalence from (\ref{f:homeq_J})).
\enpr

\Prf
The uniqueness part follows immediately from
Proposition \ref{p:noet}, since both
$
 C_*^-(f,v;\l)
$
and
$
S_*
$
are homotopy finitely generated. 
Now to  the  existence of $\mu$.
Recall that $C^-_*(f,v;\l) =
|\MM|_*$
(see (\ref{f:string_m})).
 Introduce an auxiliary string
\bq
\ZZ
=
\{S_*(1)\lTo^{p_2}
S_*(2)
{\lTo}
...
{\lTo}
S_*(n-1)
{\lTo^{p_n}}
S_*(n)
{\lTo}
...
\}
\end{equation}
where
the map
$p_n:S_*(n)\to S_*(n-1)$
is induced by the inclusion of the pairs
$
(  V^-, {t^nV^-})
\sbs
(  V^-, {t^{n-1}V^-}).
$
As the first step we shall prove 
$|\MM|_*\sim  |\ZZ|_*$.
See \cite{patou}, p.325 
for the proof of the next lemma.
\bele  \label{l:hom_com}
For every $n$ the following square is homotopy commutative
\begin{diagram}[size=2em,LaTeXeqno]
\MM_*(n) & \lTo^{\pi_{n+1}} & \MM_*(n+1)    &   \\
\dTo_{J_n} &            &\dTo_{J_{n+1}}   &\quad\lb{f:hom_com}\\
S_*(n)    & \lTo^{p_{n+1}} & S_*(n+1)& \hspace{2cm} \square\\
\end{diagram}
\enle
Therefore we have  a map  of strings
$
J:\MM\to\ZZ
$
which induces by \ref{p:map_strings}
a \heeq
~$
|J|
:
C_*^-(f,v;\l) \to |\ZZ|_*
$
of chain complexes over $\wh L_-$.
Now we are going to construct a \heeq~
$
S_*
\to |\ZZ|_*$.
The quotient chain maps
$S_*
\to
S_*(n)
=
S_*\tens{\wh L_-} L_n$
are compatible with each other, so
they induce a chain map
\bq
\xi: S_*\to\liminv~ S_*(n)
=
|\ZZ|_*
\end{equation}

\bele
The map $\xi$
is a \heeq.
\enle
\Prf
We have the following commutative diagram

\begin{diagram}[size=1.7em,LaTeXeqno]
S_* & \quad\rTo^{\xi} &       \liminv~S_*(n)        \\
\uTo &   &\uTo \quad           \\
D_* & =     &  \liminv~C_*^\D( V^-, {t^n V^-} )       \\
\end{diagram}
Both vertical arrows are \heeq s
and the lemma follows. $\qs$

Apply now Proposition
\ref{p:models}
to obtain a chain \heeq~
$\mu:
C_*^-(f,v;\l)
\to  S_*
$,
\sut
~ $\xi\circ\mu\sim |J|$.
Since both
$ C_*^-(f,v;\l) $
and
$S_*$
are free chain complexes over
$\wh L_-$,
the map $\mu$
is a chain \heq~.
The property
$\mu_n\sim J_n$
goes by construction,
and the
 existence of $\mu$ is proved.
Now the construction of the homotopy equivalence
$\Phi$ from the Theorem A is straightforward.
Recall that $\wh L=S^{-1} \wh L_-$
where $S$ is the multiplicative subset
$S=\{t^n~|~n\in \NNN\}$.
Also
\begin{gather*}
C_*(f,v)=S^{-1} C_*^-(f,v;\l)
\\
C_*^s(\bar M)\tens{ L}\wh L
=
C_*^s( V^-) \tens{ L_-}\wh L
=
S^{-1}
\Big( C_*^s( V^-) \tens{ L_-}\wh L_- \Big)
\end{gather*}
Therefore localizing $\mu$
we obtain a chain homotopy equivalence
$$
S^{-1}\mu : C_*(f,v)\rTo^\sim 
C_*^s(\bar M)\tens{ L}\wh L
$$
as required in Theorem A.
It remains to check that the chain homotopy
equivalence thus constructed depends only 
on the triple $(M,f,v)$
and not on the auxiliary choices we have made during the construction.
To this end, note first of all 
that for given $(M,f,v)$
the chain homotopy class of $S^{-1}\mu$
depends apriori only on the choice 
of the lift $F:\bar M\to\RRR$ and 
of the regular value
$\l$ (by Prop. \ref{p:mu}).
We shall therefore denote the homotopy class of $S^{-1}\mu$ by
$\Phi(F, \l)$.
Our next aim is to show that actually
$\Phi(F, \l)$
does not depend on the choice of $F$ and $\l$ neither.
It is clear that $\Phi(F-1,\l)=\Phi(F,\l+1)$, thus it
suffices to prove the next proposition.
\bepr\label{p:indep_l}
Let $\l_1, \l_2$ 
be regular values of $F$, then
$\Phi(F,\l_1)\sim \Phi(F,\l_2)$.
\enpr
\Prf
Assuming  $\l_1<\l_2$, set
$
V_{(1)}
=
F^{-1}(\l_1),\quad
V_{(2)}
=
F^{-1}(\l_2)$ and
\bq
\qquad
V_{(1)}^-
=
F^{-1}(]-\infty, \l_1]), \quad
V_{(2)}^-
=
F^{-1}(]-\infty, \l_2]).
\end{equation}
Consider the inclusions
\bq
I:
V_{(1)}^-
\rInto
V_{(2)}^-
,\qquad
I_n
:
(
V_{(1)}^-
,
t^n
V_{(1)}^-
)
\rInto
(
V_{(2)}^-
,
t^n
V_{(2)}^-
)
\end{equation}
We have also  a chain map of the corresponding Morse complexes
\bq
(I_n)_!:
\MM_*^{(1)}(n)
\to
\MM_*^{(2)}(n)
\end{equation}

\bele\label{l:dia_modk}
The following diagram is homotopy commutative:

\begin{diagram}[size=2.2em,LaTeXeqno]
\MM_*^{(1)}(n)  & \rTo^{(I_n)_!} & \MM_*^{(2)}(n)
 \\
\dTo^{J_n^{(1)}} &   &     \dTo^{J_n^{(2)}}  \quad \lb{f:dia_modk}
\\
C_*^s\Big( V_{(1)}^-,{t^n V_{(1)}^-} \Big)
& \rTo^{(I_n)_*} &
C_*^s\Big( V_{(2)}^-,{t^n V_{(2)}^-}  \Big)
\end{diagram}
where
$J_n^{(1)},
J_n^{(2)},$
are the chain homotopy equivalences from
(\ref{f:homeq_J}).

\enle

\Prf
Let
$
W_{1,n}
=
F^{-1}([\l_1-n, \l_1]),~
W_{2,n}
=
F^{-1}([\l_2-n, \l_2])
$.
The  maps $
J_n^{(1)},
J_n^{(2)}
$
are constructed with the help of $t$-ordered Morse functions
say
$g_1
,
g_2
$
on
$W_{1,n},~
W_{2,n}$.
The function $g_1$ defines a filtration in the pair
$
(
V_{(1)}^-
,
t^n
V_{(1)}^-
)
$;
let
$(    X_1^{(k)},      t^nV_{(1)}^-)
$
denote the $k$-th term of this filtration.
Similarly, we obtain a filtration
$(X_2^{(k)}
 , t^nV_{(2)}^-)
$
in
$( V^-_{(2)},
t^n V_{(2)}^-
)$.
Proposition \ref{p:thin_fi}
implies that we can choose
$g_1$
so that
$X_1^{(k)}$
be arbitrary close to
$$
t^nV_{(1)}^-
\cup
\bigcup_{\substack{p\in W_n, \\ \ind p\leq k}}
D(p;v).
$$
In particular, we can assume
that
$
X_1^{(k)}
\sbs
X_2^{(k)}     
   $.
Then the map
$I_n$
preserves the corresponding filtrations in singular homology,
and our result follows from
Lemma
\ref{l:cellular}. $\qs$

\beco\label{c:presque}
The following diagram is homotopy commutative
\begin{diagram}[size=1.8em,LaTeXeqno]
C_*^-(f,v;\l_1)  & \rTo^{I} & C_*^-(f,v;\l_2)
\\
\dTo_{\mu_1} &   &\dTo_{\mu_2}  \quad \lb{f:dia_min}
\\
C_*^s({V_{(1)}^-})\tens{ L_-}\wh L_-
& \rTo^{I_*} &
C_*^s({V_{(2)}^-})\tens{L_-}\wh L_-
\end{diagram}
(here $\mu_1, \mu_2$
are the homotopy equivalences corresponding 
to $\l_1$, resp. $\l_2$ 
by \ref{p:mu}).
\enco
\Prf
The reduction modulo
$t^n$
of this diagram coincides with
(\ref{f:dia_modk})
and, therefore, our Corollary
follows from
\ref{co:square}.
$\qs$

Now Proposition \ref{p:indep_l}
follows immediately since after localizing the square
(\ref{f:dia_min})
both horizontal arrows become identity maps, and 
localized 
vertical arrows are exactly $\Phi(\l_1)$ and 
$\Phi(\l_2)$.
$\qs$

\subsection{Functoriality of $\Phi(M,f,v)$}
\label{su:functo}

We proceed 
to the proof of the second part of Theorem A.
Let
$
F_1:\bar M_1\to\RRR,
F_2:\bar M_2\to\RRR
$
be lifts of $f_1, f_2$.
The diffeomorphism $\bar g$
sends  $S_k(F_1)$ to $S_k(F_2)$,
therefore for every $k$ it determines an isomorphism
of $\wh L$-modules $\bar g_!: C_k(f_1,v_1)\to C_k(f_2,v_2)$,
which commutes with
boundary operators, since $\bar g$
sends the orbits of $v_1$ to those of $v_2$
(by condition $\b$).
Choose a regular value $\l_1$
for $F_1$, and a regular value
$\l_2$ for $F_2$,
set
$
V_{(1)}
=
F_1^{-1}(\l_1),\quad
V_{(2)}
=
F_2^{-1}(\l_2)$ and
\bq
\qquad
V_{(1)}^-
=
F_1^{-1}(]-\infty, \l_1]), \quad
V_{(2)}
=
F_2^{-1}(]-\infty, \l_2]).
\end{equation}
Assume that
$\l_1, \l_2$
are chosen in such a way, that
\bq\label{f:cond}
\bar
g(V_{(1)}^-)
\sbs
V_{(2)}^-
\end{equation}
For $i=1,2$ let $\mu_i:
C_*^-(f_i,v_i;\l_i)\to C_*^s(V_{(i)}^-)\tens{L_-}\wh L_-$
be the corresponding chain homotopy equivalences
(see Proposition \ref{p:mu}).
The commutativity of the diagram
(\ref{f:squa_theo}) follows from the next proposition.

\bepr\label{p:presque_2}
The following diagram is homotopy commutative:
\begin{diagram}[size=1.8em,LaTeXeqno]
C_*^-(f_1,v_1;\l_1)  & \rTo^{\bar g_!} & C_*^-(f_2,v_2;\l_2)
\\
\dTo_{\mu_1} &   &\dTo_{\mu_2}  \quad \lb{f:dia_min_2}
\\
C_*^s({V_{(1)}^-})  \tens{L_-}\wh L_-
& \rTo^{\bar g_*} &
C_*^s({V_{(2)}^-})  \tens{L_-}\wh L_-
\end{diagram}
\enpr

 The proof of this proposition is completely similar to that
of Corollary \ref{c:presque}
and will be omitted. 
Now, applying the $S$-localization to the square (\ref{f:dia_min_2})
we obtain the second part of Theorem A. $\qs$

\section{$C^0$-perturbations of gradient-like vector fields}
\label{s:cond_c}

This section contains recollections
of results from \cite{pafest}, which are
crucial for the proof of Theorem B.
In Subsections 
\ref{su:delta_thin} -- \ref{su:c_hand}
we work with real-valued functions of cobordisms
and their gradients.
The first and second subsections contain 
the results about 
behaviour of stable manifolds and similar objects under 
$C^0$-small perturbations of the gradient.
In Subsection \ref{su:cond_c}
we recall the condition $(\gC)$
on the gradient-like vector field 
and the properties of 
\glvf s satisfying this condition.
Subsections
\ref{su:grad_desc} --
\ref{su:c_hand}
contain an account of the properties of $f$-gradients
satisfying $(\gC)$.
In 
Subsection \ref{su:fil_cov}
we proceed to the  circle-valued Morse
maps and we 
recall the condition $(\gC')$
which is the appropriate analog of $(\gC)$
in this framework.
In Subsection \ref{su:psi}
we give the construction of the chain equivalence $\psi$
from \cite{pafest}
between the Novikov complex and
the completed singular chain complex of $\bar M$.

\subsection{$\d$-thin handle decompositions}
\lb{su:delta_thin}

Through  Subsections \ref{su:delta_thin} -- \ref{su:c_hand}
$\fcob$ is a Morse function on a riemannian
cobordism $W$ of dimension $m$, and
$v$ is an $f$-gradient.
We shall need some more of Morse-theoretic terminology.
Set 
$$U_1=\{x\in \pr_1 W~|~\g(x,\cdot;v) \mbox{ reaches } \pr_0 W\}.
$$
Then $U_1$ is an open subset of $\pr_1 W $
and the gradient descent along the trajectories of $v$
determines a diffeomorphism of $U_1$
onto an open subset $U_0$
of $\pr_0 W$.
This diffeomorphism will be denoted by
$\stv$
and we shall abbreviate
$\stv (X\cap U_1)$ to
$\stv (X)$.
We say that
$v$ satisfies {\it Almost Transversality Condition}
if
   $$\big( x,y \in S(f) ~\&~
\ind x \leq \ind y \big) \Rightarrow \big(
D(x,v) \pitchfork D(y,-v) \big)$$
\vskip0.1in
 A Morse function $\phi:W\to[\a,\beta]$
is called {\it adjusted to $(f,v)$}, if:

 1) $S(\phi)=S(f)$, and $v$ is also a $\phi$-gradient.

 2) The function $f-\phi$ is constant in a \nei~
of $\pr_0W$, in a \nei~ of $\pr_1 W$, and in a \nei~ of each
point of
$S(f)$.

Using the standard techniques of rearrangement of critical points
(see for example \cite{milnhcob}, \S 4, or an 
exposition which uses our present terminology in \cite{paadv}, \S 2)
it is not difficult to
show that for an arbitrary Morse function $f$
and an $f$-gradient satisfying the the~\ata,
there is an ordered Morse function $g$, adjusted to $(f,v)$.
Let $p\in W\sm \pr W$, and  $\d>0$.
Assume that for some $\d_0>\d$ the restriction of
the exponential map
$\exp_p:T_pW\to W$ to the disc $B^m(0,\d_0)$
is a \dfm~on its image.
Denote by $B_\d(p)$ (resp. $D_\d(p)$) the riemannian open ball
(resp. closed ball)
of radius $\d$
centered in $p$. We shall use the notation
$B_\d(p), D_\d(p)$ only when the assumption above on $\d$ holds.
Set
\begin{gather*}
B_\delta(p,v)=\{x\in W~\vert~ \exists
t\geq 0 :
 \gamma (x,t;v)\in  B_\delta (p)\}\\
D_\delta(p,v)=\{x\in W~\vert ~\exists
t\geq 0 :
\gamma (x,t;v)\in  D_\delta (p)\}
\end{gather*}

We denote by $B_\d (\indl s;v)$, 
the union of
$B_\d(p,v)$, where $p$ ranges over critical
points of $f$ of index $\leq s$.
We shall also use similar
notation like $B_\d (\indl s;v)$, $D_\d(\inde s;v)$ or
$B_\d(\indg s;v)$, which is now clear
without special definition.
(Actually the symbol 
$D(\indl s;v)$
was already introduced
in 
\rref{f:indl_s}.)
 Set
$$C_\d(\indl s; v)=W\sm B_\d(\indl {m-s-1};-v)$$

Let $\phi:W\to[a,b]$ be an ordered Morse function with an
ordering sequence
$a_0<a_1<\dots <a_{m+1}$. Let $w$ be a $\phi$-gradient.
Denote $\phi^{-1}\([a_i,a_{i+1}])$ by $W_i$.

\begin{defi}
We say that $w$ is {\it $\d$-separated \wrt~
$\phi$ }
(and the ordering sequence $(a_0,\dots,a_{m+1})$), if

i) for every $i$ and every $p\in S_i(f)$ we have
$D_\d(p)\sbs W_i\sm \pr W_i$;

ii) for every $i$ and every $p\in S_i(f)$
there is a Morse function
$\psi:W_i\to[a_i,a_{i+1}]$, adjusted to
$(\phi\mid W_i, w)$
and a regular value $\l$ of $\psi$ \sut~
$$D_\d(p)\sbs\psi^{-1}(]a_i,\l[)$$
and for every $ q\in S_i(f), q\not= p$ we have
$$D_\d(q)\sbs \psi^{-1}(]\l,a_{i+1}[)$$

\end{defi}

We say that $w$ is {\it $\d$-separated }if
it is $\d$-separated \wrt~some ordered Morse
function $\phi:W\to [a,b]$, adjusted to $(f,v)$.
Each $f$-gradient satisfying \ata~ is $\d$-separated
for
some $\d>0$.

\bepr [\cite{pafest}, Prop. 3.2, 4.1]
\lb{p:ddd}

If $v$ is  $\e$-separated, then
$\forall
\delta\in[0,\e]$ and
$\forall
s: 0\leq s\leq m$
\been
\item $D_\delta (\indl s;v)$         is compact.
\item 
The family 
${\bigcap} _{\theta>\d} B_\theta (\indl s;v)
$
form a fundamental system of \nei s
of 
 $D_\delta (\indl s;v)$.
The family 
${\bigcap} _{\theta>0} B_\theta (\indl s;v)
$
form a fundamental system of \nei s
of 
 $D (\indl s;v)$.

\item
 $B_\delta (\indl s;v)
=\Int D_\d(\indl s; v)$
and
$ D_\d(\indl s; v)
\sbs
 C_\d(\indl s; v)$.
\item
The group
$$
H_*\Big(
 D_\d(\indl s; v)\cup \pr_0 W,~
 D_\d(\indl {s-1}; v)\cup \pr_0 W
\Big)
$$
vanishes  if $*\not= s$
and is a free abelian group generated by the
classes of the descending discs $D(p,v)$
with $p\in S_s(f)$.
\end{enumerate}
\enpr

We shall often denote
$D_\d(\indl s;v)$
by $W^{\{\leq s\}}$
if the values of $v,f,\d$
are clear from the context.

\subsection{$C^0$-continuity properties of descending discs}
\label{su:c0_cont}
In this subsection we study the behaviour of descending discs
under $C^0$-small perturbations of the gradient.
let $G(f)$
denote the set of all $f$-gradients.

\bele\label{l:cont_1}
Let $\d>0$
and $K\sbs B_\d(v)$
a compact set. There exists $\eta>0$, \sut~
for every $w\in G(f)$
with $||w-v||<\eta$
we have $K\sbs B_\d(w)$.
\footnote 
{In the present paper  the symbol
$||\cdot||$ denotes the $C^0$-norm.}
\enle
\Prf
Set $B=\cup_{p\in S(f)} B_\d(p)$.
For a subset $Q\sbs W$
let us denote by
$\RR(Q, w)$
the following condition:
\begin{quote}
For every $x\in Q$
the trajectory
$\g(x,t;w)$
intersects $B$.
\end{quote}
Note that $\RR(K,v)$ holds. For
$x\in K$
choose $t(x)\geq 0$
\sut~
$\g(x,t(x);v)\in B$.
By Corollary 5.5 of \cite{pastpet}
there is a \nei~$S(x)$
of $x$ and $\d(x)>0$
\sut~ for every
$f$-gradient $w$ with $||w-v||<\d(x)$
the property $\RR(S(x),w)$
holds. Choosing a finite covering of $K$
by the subsets $S(x)$ finishes the proof. $\qs$

\bele\label{l:cont_2}
Let $\d>0$
and assume that $v$ is $\d$-separated.
Let $U\sbs W$
be an open set and assume 
$D_\d(\indl k;v)\sbs U$.
Then there is $\e>0$
\sut~
for every $w\in G(f)$ with $||w-v||<\e$
we have
$D_\d(\indl k;v)\sbs U$.
\enle
\Prf
First we prove the lemma 
for $k=\dim M$,
so that
$D_\d(\indl k;v)$
is the union of all $\d$-thickened descending discs.
Let $K=W\sm U$.
Let us denote by $Q(w)$ the following condition:
\begin{quote}
For every $x\in K$ the trajectory $\g(x,t;w)$
reaches $\pr_1 W$,
without intersecting the set
$\bigcup_{p\in S(f)} D_\d(p)$.
\end{quote}
We know that $Q(v)$ holds
and we have to prove that $Q(w)$
holds for every \glvf~$w$
sufficiently $C^0$-close to $v$.
This follows directly from Corollary
5.6 of \cite{pastpet}.
To prove the lemma for
any $k$ we consider an ordered Morse function
$\phi:W\to[a,b]$
adjusted to $(f,v)$
and we set
$W'=\phi^{-1}([a_0,a_{k+1}])$.
Then
$D_\d(\indl k;v)
=
D_\d(\indl {m}; v|W')$
and applying the previous argument to $W'$ we are over. $\qs$

\subsection{Condition $(\gC)$}
\lb{su:cond_c}

In this subsection we recall the condition
$(\gC)$
on the gradient $v$.
If this condition holds, the gradient descent map, corresponding
to
$v$ can be endowed with a structure, resembling closely
the cellular maps between $CW$-complexes.
Also the cobordism $W$ can be endowed with
certain handle-like filtration
which encompasses the handle-like filtrations on $\pr_0 W$
and $\pr_1 W$.

\bede\lb{d:cc}
(\cite{pafest}, Def. 4.5)
We say that
{\it $v$ satisfies condition $(\gC)$}
if there are objects 1) - 4), listed below,
with the properties (1 - 3) below.
\pa
{\it Objects:}
\pa
\begin{enumerate}
\item[1)] An ordered Morse function $\phi_1$
on $\pr_1 W$ and a $\phi_1$-gradient $u_1$.
\item[2)]  An ordered Morse function $\phi_0$
on $\pr_0 W$ and a $\phi_0$-gradient $u_0$.
\item[3)]  An ordered Morse function $\phi$
on $ W$
adjusted to $(f,v)$.
\item[4)]  A number $\d>0$.
\end{enumerate}
\pa
{\it Properties:}
\pa
(1)\quad  $u_0$ is $\d$-separated \wrt~ $\phi_0$,
$u_1$ is $\d$-separated \wrt~ $\phi_1$, $v$ is
$\d$-separated \wrt~ $\phi$.
\pa

 \begin{multline*}(2)\qquad    \stv \bigg(C_\d(\indl j;u_1)\bigg)
\cup
\bigg(D_\d\(\indl {j+1};v\)\cap \pr_0 W\bigg)
\sbs                \\
B_\d(\indl j,u_0)
\mbox{ for every } j
\end{multline*}

\begin{multline*}(3) \qquad     \st v\bigg(C_\d(\indl j;-u_0)
\bigg)
\cup
\bigg( D_\d(\indl {j+1}; -v)\cap \pr_1 W\bigg)
\sbs                  \\
B_\d(\indl j;-u_1)
\mbox{ for every } j
\end{multline*}
$\qt$
\end{defi}

The set of all $f$ gradients satisfying $(\gC)$ will be denoted by
$GC(f)$.

\pa

\begin{theo}\lb{t:cc}
(\cite{pafest}, Th. 4.6)
$GC(f)$ is open and dense in $G(f)$ \wrt~  $C^0$ topology.
Moreover, 
let $v_0\in G(f)$, let $U$ be a \nei~of $\pr W$, and $\d>0$.
Then there is $v\in GC(f)$ with
$||v-v_0||<\d$
and
$\supp(v-v_0)\sbs U$.
\end{theo}

\subs{Condition $(\gC)$ and
regularization of gradient descent map
       }
\label{su:grad_desc}

As we have already mentioned the application
$\stv$
is not everywhere defined. But if the gradient $v$
satisfies the condition
$(\gC)$,
we can associate to $v$
some family of continuous maps which plays the role of
"cellular approximation" of $\stv$,
and a homomorphism
({\it homological gradient descent})
which is a substitute for the
homomorphism induced by $\stv$ in homology.
Let $v$ be an $f$-gradient satisfying $(\gC)$.
We have two ordered Morse functions
$\phi_1, \phi_0$
on
$\pr_1 W$, resp.  $\pr_0 W$,
and their gradients $u_1, u_0$,
which give rise to
filtrations on $\pr_1 W$, resp. $\pr_0 W$.
Namely,
let
$\b_0< \b_1< ... <\b_m$
be an ordering sequence for
$\phi_1$,
 $\a_0< \a_1< ... <\a_m$
be an ordering sequence for
$\phi_0$.
Let
\bq\label{f:defidefi}
\dow^{(k)}
=
\phi_1^{-1}([\b_0,\b_{k+1}]), \quad
\daw^{(k)}
=
\phi_0^{-1}([\a_0,\a_{k+1}])
\end{equation}
For $k\geq 0$ consider the set $Q_k$
of all
$x\in
\daw^{(k)}
$
where $\stv$ is {\it not} defined.
Equivalently,
$$
Q_k=\{x\in \daw^{(k)}~
|~\g(x,t;-v) \mbox{ converge  to point in } S(f)
\mbox{ as } t\to\infty\}.
$$
This is a compact set, and the  condition $(\gC)$
implies that
 $Q_k\sbs D(\indg {m-k} ; -v)$.
Therefore there is a \nei~ $U$
of $Q_k$ in
$
\daw^{(k)}
$
\sut~
$\stv (U)$
is in
$D_\d(\indl k; v)\cap \dow$
and this last set is in
$\Int \daw^{(k-1)}
$
(again by $(\gC)$).
It follows that the map
$\stv$
gives rise to a well-defined continuous map
\bq
\vflesh:
\daw^{(k)}
/
\daw^{(k-1)}
\to 
\dow^{(k)}
/
\dow^{(k-1)}
\end{equation}
The homomorhism, induced by this map in homology,
is called {\it homological gradient descent}.
It
will be denoted by
$\HH_k(-v)$
or simply $\HH_k$
if the gradient $v$ is clear from the context.

\subsection{Morse type decomposition of
 cobordisms: preliminary discussion}
\label{su:morsedec_prel}

For any Morse function on a closed manifold
the stable manifolds of the critical points 
(\wrt~to a gradient-like vector
field satisfying the transversality assumption)
form a cellular-like decomposition of the manifold.
This construction
 generalizes directly to Morse functions on  cobordisms,
but we lose some of its properties. Namely, 
the union of descending discs of critical  points
is not necessarily equal to $W$
(consider for example the second coordinate projection
$V\times [0,1]\to [0,1]$;
this is a Morse function without critical points).
In this subsection we discuss a natural approach 
to overcome this difficulty.

\bede\label{d:track}
For $X\sbs \pr_1 W$
let us denote
by $T(X,-v)$ and call {\it track of $X$}
the
union of all
$(-v)$-trajectories
starting at a point of the set $X$.
Similarly, for
$Y\sbs \dow$
let
$T(Y,v)$
be
union of all
$v$-trajectories
starting at a point of the set $Y$.
\end{defi}

A nice cellular-like decomposition
fitted to $f$ and $v$ should contain:\label{cell_like}
\been\item The cells of $\pr_0 W$,
\item the cells of $\pr_1 W$,
\item
the descending discs $D(p,v)$ for $p\in S(f)$,
\item the tracks of the cells of $\pr_1 W$.
\enen
By "cells of $\pr_1 W$"
we mean descending discs corresponding to some
Morse function
$\phi_1:\pr_1 W\to\RRR$
and its gradient $u_1$
(similarly for $\pr_0 W$).
In order to  get a cellular-like  decomposition,
some requirements must be met.
Namely, for $p\in S(f)$
the sole
$D(p,v)\cap\pr_0 W$
of
the descending disc $D(p,v)$
with $p\in S(f)$
must belong to the
$(\ind p -1)$-skeleton
of $\pr_0 W$.
Let us write down this condition:
\bq\label{f:requir}
D(p,v)\cap\pr_0 W
\sbs
\pr_0 W^{(\ind p -1)}
\end{equation}
Similar requirement must hold for the tracks of the cells of $\pr_1 W$.
Although it is plausible that such requirement could be met
(see the discussion in \cite{pafest}, end of Section 2)
I do not know whether it is true in general.
The method
which we apply to realize this idea 
uses the thickenings of descending discs.
This is how the conditions of the type $(\gC)$
appeared first, in  \cite{pastpet}.
 We introduced in this paper (\S 4)
the conditions $(RP)$
on the $f$-gradient, 
which  allow to regularize the gradient descent map
and to give it a cellular-like structure.
The condition $(\gC)$, introduced
in \cite{pafest}, being quite similar to $(RP)$,
has some technical advantages over it, for example 
$(\gC)$ is  an open and dense condition \wrt~ $C^0$-topology
(while $(RP)$ was only proved to be $C^0$-dense).

If we look now  at the condition
$(\gC)$ from the point of view of 
\rref{f:requir}, we notice that the part $(2)$
of $(\gC)$
is exactly  the "$\d$-thickened" analog of this requirement.
In the next subsection we shall see that the condition
$(\gC)$ gives rise indeed to a filtration
of the cobordism $W$, 
which has the properties of handle
filtrations of closed manifolds.

Our main aim is in the applications of these methods to
the circle-valued maps. To this end we note that 
in the framework of
circle-valued Morse maps the corresponding analog of
the  requirement 
\rref{f:requir} can not be met, while 
the methods based on the condition $(\gC)$ 
work also for this generalization. These methods 
will be our main tools
for the construction of the chain
equivalence from Theorem B.

\subsection{Handle-like filtrations of cobordisms}
\label{su:c_hand}

In this subsection we associate to
 each $f$-gradient $v$ satisfying the condition $(\gC)$
a handle-like filtration
of the cobordism.
Let $\phi$
be the ordered Morse function on $W$, from the definition
\ref{d:cc}, and
$(a_0,..., a_{m+1})$
be an ordering sequence for $\phi$.

Let $Z_k$
be the set of all
$z\in
 \phi^{-1}([a_0,a_{k+1}])$
\sut
~the $(-v)$-trajectory
$\g(z,t;-v)$
converges to a critical point
of $\phi$
as $t\to\infty$
in
$ \phi^{-1}([a_0,a_{k+1}])$
or reaches $\pr_0 W$
and intersects it at a point
in
$\pr_0 W^{(k)}$.
In other words
\bq
Z_k=\Big(T(\dow^{(k)}, v) \cup D(-v)\Big)
\cap
\phi^{-1}([a_0, a_{k+1}])
\end{equation}
(where $D(-v)$
stands for the union of all ascending discs of $v$).
Note in particular that $Z_k$
contains the set
$\pr_0 W^{(k)}$
and the descending discs
$D(p,v)$
for $\ind p\leq k$.

Let $T_k$
be the set of all points
$y\in \phi^{-1}([a_{k}, b])$
\sut~the $v$-trajectory
starting at $y$ reaches
$\pr_1 W$
and intersects it at
$\pr_1 W^{(k-1)}$.
In other words
\bq
T_k
=
T(
\pr_1 W^{(k-1)}; -v)
\cap
\phi^{-1}([a_{k}, b])
\end{equation}
Now we can define the filtration
$\{W^{\langle k \rangle}\}$
of $W$; set
\bq\label{f:fil_c}
W^{\langle k \rangle}
=
\pr_1 W^{(k)}\cup
Z_k
\cup
T_k.
\end{equation}
Thus
$
W^{\langle k \rangle}
$
contain the handles
of index $\leq k$
of the manifolds
$\pr_0 W, ~\pr_1 W,~ W$,
and also the tracks of the handles of index $\leq k-1$
of $\pr_1 W$.
(The sets $Z_k, T_k, W^{\langle k \rangle}
,
W^{\langle k-1 \rangle}
$ are depicted  on the page \pageref{figure}.)
Consider the corresponding filtration in
the singular homology and
denote by $E_*$ the corresponding adjoint complex, so that
$$
E_k = H_k
(
W^{\langle k \rangle}
,
W^{\langle k-1 \rangle}
).
$$
In \cite{pafest}
we proved that the filtration induced by $
W^{\langle k \rangle}$
in
the singular
homology  is nice, and   computed the  boundary operator in $E_*$
 in terms of Morse complexes of
$\phi_0, \phi_1$ and $\phi$
and the gradient descent homomorphism.
We shall now recall these results.
Consider the Morse complexes
$C_*(\phi_0, u_0),~
C_*(\phi_1, u_1),~
C_*(\phi, u)$
associated to Morse functions $\phi_0, \phi_1, \phi$.
The obvious inclusions

\bq
\Big(
\dow^{(k)}
,
\dow^{(k-1)}
\Big)
\rInto
(
W^{\langle k \rangle}
,
W^{\langle k -1\rangle}
)
\lInto
\Big(
\daw^{(k)}
,
\daw^{(k-1)}
\Big)
\end{equation}
induce chain maps
\bq
C_*(\phi_0,u_0)
\rTo^{\l_0}
E_*
\lTo^{\l_1}
C_*(\phi_1,u_1)
\end{equation}
It follows from the
condition $(\gC)$
that for every $p\in S_k( f)$
the descending disc
$D(p,v)$
has a well defined fundamental class $[p]$
in
$
H_k
(
\wmok
,
\wmokm
)
$
and the corresponding map
$p\mapsto [p]$
defines an inclusion
\bq
C_k(\phi,v)\rTo^{\mu}
H_k
(
\wmok
,
\wmokm
)
\end{equation}
The images of $\l_1, \l_0$ and $\mu$
correspond  to the components 1) -- 3)
of the hypothetical cell-like decomposition from the page
\pageref{cell_like}.
Now to the fourth
component.
The condition $(\gC)$
implies that no $(-v)$-trajectory starting at
$x\in
\pr_1 W^{(k-1)}$
converges to a critical point of $\phi$
in
$\phi^{-1}([a_{k}, b])$,
therefore we have a homeomorphism
\bq
T_k\approx
\pr_1 W^{(k-1)}
\times
I
\end{equation}
where $I=[a_{k},b]$,
and
the pair
$
(T_k,
T_k\cap
W^{\langle k -1 \rangle}
)
$
is therefore homeomorphic
to
\bq
\Big(
\dwmokm, \dwmokmm)
\Big)
\times
(I,\pr I).
\end{equation}
The multiplication with
the fundamental class of
$(I, \pr I)$
defines the map
\bqq
\tau
:
H_{k-1}\Big( {\dwmokm}, {\dwmokmm}\Big)
=
C_{k-1}(\phi_1, u_1)
\to H_k({\wmok}, {\wmokm})
=
E_k
\end{equation*}

\bepr[\cite{pafest}, Theorem 5.5]\label{p:hom_grw}
$H_*(\wmok, \wmokm)=0$
for $*\not=k$, and the map
\begin{multline}
\label{f:decomp_w}
\LL=(\l_1, \l_0, \mu, \tau)
:   \\
C_k(\phi_1, u_1)\oplus
C_k(\phi_0,u_0)\oplus C_k(\phi,v)
\oplus
C_{k-1}(\phi_1, u_1)  \rTo \\
\rTo H_k( \wmok, \wmokm)
\end{multline}
is an isomorphism.
\enpr

\bepr[\cite{pafest}, Proposition 5.9]
The matrix of the  boundary
operator
$\pr: E_k\to E_{k-1}$
 \wrt~the decomposition
(\ref{f:decomp_w})
is

\bq\label{f:bou_w}
\left(
\begin{matrix}
 \pr_{k}^{(1)}   & 0 & 0 & 1 \\
0 &  \pr_{k}^{(0)} & * & -\HH_{k-1}    \\
0 & 0 &  \pr_{k}   & * \\
0 & 0 & 0 &  -\pr_{k-1}^{(1)} \\
\end{matrix}   \right)
\end{equation}

\enpr
Here
$\HH_{k-1}$ is the homological gradient descent
homomorphism (Subsection \ref{su:grad_desc}), and 
$
 \pr_{*}^{(1)},~~  \pr_{*}^{(0)},
~ \pr_{*}
$
are the boundary operators
 in the Morse complexes
$
C_*(\phi_1, u_1)
$, resp. $
C_*(\phi_0, u_0),
$
and
$
C_*(\phi, u)
$.
 The terms denoted by  $*$
are not important for us.

\subsection{Filtration of $\bar M$
associated with $C^0$-generic gradient}
\label{su:fil_cov}

We proceed now to circle-valued Morse functions.
Let $f:M\to S^1$
be a Morse map,
$v$ be an $f$-gradient.
We shall work with the terminology of Section
\ref{su:morse_wn}.
Assume that $M$ is endowed with a riemannian metric;
then
$\bar M$
and
$W=F^{-1}([\l-1,\l])$
inherit the riemannian structure.
The generator $t$ of the structure group
of the covering
$\bar M\to M$
is an isometry of $\bar M$
and induces an isometry
$\pr_1 W\to \pr_0 W$.

We shall say that $v$ satisfies condition
$(\gC')$
if the $(F\vert W)$-gradient $v$ satisfies the condition $(\gC)$
from \su~\ref{su:cond_c},
and, moreover, the Morse functions $\phi_0, \phi_1$
and their gradients
$u_0, u_1$
can be chosen so as to
satisfy
$\phi_0(tx)=\phi_1(x), t_*(u_1)=u_0$.
The set of
Kupka-Smale
$f$-gradients $v$ satisfying $(\gC')$
will be denoted by
$\GG_0(f)$.
The set $\GG_0(f)$
is $C^0$-open and dense in $\GG(f)$
(this is a version of the theorem
\ref{t:cc}, see
\cite{pafest}, \S 8).
Let $v\in \GG_0(f)$.
Consider the corresponding filtration
$
\{W^{\langle k  \rangle}\}
$
of $W=F^{-1}([\l-1,\l])$
from the section \ref{su:c_hand}.
Define the filtration
$
\{V^-_{\langle k  \rangle}\}
$
of
$V^-$
by
\bq
V^-_{\langle k  \rangle}
=
\bigcup_{s\geq 0} t^s
W^{\langle k  \rangle}
\end{equation}
We shall now  describe the homology
$
H_*
(
 V^-_{\langle k  \rangle}
,
 V^-_{\langle k-1  \rangle}
)
$.
Note first of all the splitting
\bq
\label{f:split}
  V^-_{\langle k  \rangle}
/
  V^-_{\langle k-1  \rangle}
\approx
\bigvee_{s\geq 0}
t^s
\Big(
 W^{\langle k  \rangle}
/
 W^{\langle k-1  \rangle}
\Big)
\end{equation}
It follows immediately from
\ref{p:hom_grw}
and
(\ref{f:decomp_w})
that the filtration
 $
 V^-_{\langle k  \rangle}
$
of
$
 V^-
$
is nice.
The homology
$
H_k(
  V^-_{\langle k  \rangle}
,
  V^-_{\langle k-1  \rangle}
)
$
can
be described explicitly.
 Set
\bq
\RR_*
=
C_*(\phi_1, u_1)\tens{\ZZZ } L_-, \quad
\NN_*
=
C_*(\phi, v)\tens{\ZZZ } L_-
\end{equation}
Then $\RR_*$ and $\NN_*$ are free finitely generated
$\L_-$-complexes.
Tensoring with $ L_-$
the maps
$\l_1, \mu, \tau$
we obtain \ho s
\begin{multline}
\wi\l_1:\RR_k\to
H_k(
  V^-_{\langle k  \rangle}
,
  V^-_{\langle k-1  \rangle}
), \quad
\wi\mu:
\NN_k
\to
H_k(
  V^-_{\langle k  \rangle}
,
  V^-_{\langle k-1  \rangle}
)
,\\
\wi\tau:\RR_{k-1}
\to
H_k(
  V^-_{\langle k  \rangle}
,
  V^-_{\langle k-1  \rangle}
)
\end{multline}
The module
$C_k(\phi_0, u_0)$
is identified with
$tC_k(\phi_1, u_1)$
and therefore the homological gradient descent
$\HH_k$
can be considered after
tensoring by $ L_-$, as a \ho~
$\RR_k\to \RR_k$.
Note that its image is in $t\RR_*$.
The following two propositions are proved in
\cite{pafest}, Prop. 7.4.

\bepr\label{homgr_v}
\been\item
$H_*
(
  V^-_{\langle k  \rangle}
,
  V^-_{\langle k-1  \rangle}
)
=0
$ if $*\not= k$
\item
The map
\bq\label{f:split_v}
L=(\wi\l_1, \wi\mu, \wi\tau):
\RR_k\oplus \NN_k\oplus \RR_{k-1}
\to
H_k(
  V^-_{\langle k  \rangle}
,
  V^-_{\langle k-1  \rangle}
)
\end{equation}
is an isomorphism.         \enen
\enpr
Thus the filtration induced by
$\{ V^-_{\langle k \rangle}\}$
in the singular chain complex of $V^-$
is nice; the corresponding adjoint chain complex
will be denoted by
$\EE_*$.
\bepr
\label{p:bou_v}
The matrix of the boundary operator
$\pr:\EE_k\to \EE_{k-1}$
\wrt~the direct sum decomposition
\rref{f:split_v} is
\bq\label{f:bou_v}
\left(
\begin{matrix}
 \pr_{k}^{(1)}   & * &  1-\HH_{k-1}  \\
0 &   \pr_{k}   & * \\
0 & 0 &  -\pr_{k-1}^{(1)} \\
\end{matrix}   \right)
\end{equation}    \enpr

{}$\qs$

The direct sum decomposition \rref{f:split_v}
gives rize to a natural free base in $\EE_*$, namely 
the family of the free generators of $\EE_k$ is
is
$$
S_k(\phi_1)
\sqcup
S_k(\phi)
\sqcup
S_{k-1}(\phi_1).
$$
Corollary \ref{c:heeq}
implies the existence of a chain equivalence
$$
\EE_*\to C_*^s(V^-).
$$
Composing it with the inverse of a natural chain equivalence
$$
\iota: C_*^\D(V^-)
\to
C_*^s(V^-)
$$
(where $\D$ is the triangulation of $M$ chosen so that $V$ be a simplicial subcomplex, see
the beginning of Subsection \ref{su:equi_fi})
we obtain a chain equivalence
\bq
\label{f:composit}
(\iota)^{-1}\circ\rho:
\EE_*
\to
C_*^\D(V^-)
\end{equation}
of finitely generated free based $ L_-$-complexes.
Note, that the torsion of 
this chain equivalence vanishes, since it 
belongs to
the group
$\ove{K}_1(\ZZZ[t])$
and this group is trivial by \cite{bhs}.
\label{r:tors_triv}

\subsection{ Homological gradient descent and
zeta function}
\lb{su:h_eta}

As in the preceding subsection we assume here that
$v$ satisfies the condition $(\gC')$.
We assume 
 moreover that $v\in \GG(f)$
so that the eta and zeta functions are defined.
We have
\bepr[\cite{pafest}, \S 8]
\footnote{ In \cite{pafest} there is a sign error in this formula
corrected in \cite{pawitt}.}
\lb{p:zeta_and_h}  
\bq
\lb{f:zeta_and_h}
\z_L(-v)
=
\prod_{i=0}^m \Big(\det(1-\HH_m)\Big)^{(-1)^{m+1}}
\end{equation}
\enpr
Without reproducing the full proof of this proposition
we wish to indicate here  the main idea. 
Let $\g$ be a closed orbit of $(-v)$,
and $x$ be the intersection
of $\g$
with
$V=\pr_1 W$.
There is a unique $k$ satisfying 
$x\in 
V^{\langle k \rangle}
\sm
V^{\langle k-1 \rangle}.
$
It follows from the condition $(\gC)$
that every time the trajectory $\g$
intersects $V$, the point of intersection
belongs again to
$x\in 
V^{\langle k \rangle}
\sm
V^{\langle k-1 \rangle}.
$
Therefore the set of closed orbits of 
$(-v)$
falls into the disjoint union of subsets
$\FF_k$,
indexed by $k, 0\leq k\leq m-1$
and the subset $\FF_k$
is identified with the set of periodic orbits 
of the gradient descent map
\bq
\vflesh:
\daw^{(k)}
/
\daw^{(k-1)}
\to 
\dow^{(k)}
/
\dow^{(k-1)}
\approx
\daw^{(k)}
/
\daw^{(k-1)}
\end{equation}
Now recall that $\HH_{k}$
stands for the homomorphsim induced by 
$\vflesh$
in homology, and apply  the 
Lefschetz fixed point formula
to finish the proof.

\subsection{Homotopy equivalence $\psi$}
\label{su:psi}

Now everything is ready to describe
the chain equivalence
$\psi:C_*(f,v)\to C_*^s(\bar M)\tens{L}\wh L$.
As it was the case for the equivalence
$\Phi$ from Theorem A,
we construct first an equivalence
of chain complexes
$C_*^-(f,v;\l)$ and
$C_*(V^-)\tens{L_-}\wh L_-$.
The formula
\rref{f:xi_def} defining this chain equivalence, requires
some preliminary computations
with descending discs and their homology classes.
Let $p\in S_k(F)\cap W$.
For $r\geq 1$
set
$$\Sigma_r(p)
=
D(p,v)\cap t^r\cdot V,$$
so that  $\Sigma_r(p)$
is an oriented ($k-1$)-dimensional \sma~of $t^r \cdot V$.

\bepr\label{p:sole}
\been\item
$D(p,v)\sbs \Int
 V^-_{\langle k  \rangle}$.
\item
For every
$r\geq 1$
we have
\bq\label{f:sole}
\Sigma_r(p)
\sbs t^r \cdot \Int
\daw^{(k-1)}.
\end{equation}
\item
The set
$\Sigma_r(p)
\sm t^r \cdot \Int  \daw^{(k-2)}$
is compact.
\enen
\enpr
\Prf
From the very definition of the set
$Z_k$
we have
$D(p,v)\cap W\sbs Z_k$.
The condition
$(\gC')$
implies also
$$
\Sigma_1(p)
\sbs
\Int
\dow^{(k-1)}
=
t\cdot
\Int
\daw^{(k-1)}
.$$

Thus (\ref{f:sole})
is proved for $r=1$.
It is now easy to prove
the points 2) and 3)
by induction in $r$, using the condition $(\gC' )$,
and deduce from it
the inclusion
\bq
D(p,v)\cap t^r \cdot  W\sbs t^r \cdot  \Int
 V^-_{\langle k  \rangle}. \qquad \square
\end{equation}

It follows that
the fundamental class of
$\Sigma_r(p)$
is well defined
modulo
$ t^r \cdot  \daw^{(k-2)}$. This element will be denoted
\bqq
\s_r(p)
\in
H_{k-1}\Big(
 {t^r \cdot \daw^{(k-1)}},~
 {t^r \cdot \daw^{(k-2)}}\Big)
\end{equation*}

The next formula is easy to prove  by induction in $r$:
\bq\label{f:sigma_r}
\s_r(p)=\HH_{k-1}^{r-1}(\s_1(p)).
\end{equation}

Now let $n\in \NNN$ and
consider the chain complex
$\EE_*\tens{L_-}  L_n
=
\EE_*/ t^n\EE_* $.

Using \rref{f:split}
it is not difficult to obtain a natural isomorphism:
\bq
\EE_k/  t^n\EE_k
\approx
H_k(
 V^-_{\langle k  \rangle}
,
  V^-_{\langle k -1 \rangle}
\cup t^n V^-_{\langle k  \rangle}
)
\approx
H_k(
 V^-_{\langle k  \rangle}
\cup t^n V^-
,
  V^-_{\langle k -1 \rangle}
\cup t^n V^-
)
\end{equation}

The critical points of indices $\leq k-1$
and their descending discs are contained in
$\Int
 V^-_{\langle k -1 \rangle}
$
by \ref{p:sole},
therefore the set
$
D(p,v)\sm
\Int(
 V^-_{\langle k-1  \rangle}
\cup
t^nV^-)
$
is compact.
This implies in particular, that the
chosen orientation of
$D(p,v)$
defines a homology class
\bq
\D_p\in
H_k(
 V^-_{\langle k  \rangle}
\cup
t^n V^-
,
 V^-_{\langle k-1  \rangle}
\cup t^n  V^-
)
=
\EE_k/t^n \EE_k.
\end{equation}
Our next aim is to compute the components  of this
element \wrt~ the direct sum decomposition
(\ref{f:split_v}).

\bepr\label{p:compu_d}
Let $p\in S_k(F)\cap W$.
Then
\bq
\D_p
=
[p]
-
\sum_{r=1}^{n-1}
\wi\tau(\HH^{r-1}_{k-1}(\s_1(p)))
\end{equation}
\enpr

\Prf
Let
$\ZZ=
\cup_{0\leq s\leq n-1} t^s      W^{\langle k-1  \rangle}$.
The
set
$Q=D(p,v)\sm\Int\ZZ$
falls into the disjoint
union:
$$
Q=\sqcup_i Q_i \quad \mx{ with } \quad 
Q_i
=
Q\cap
  F^{-1}([\l-i-1, \l-i])
$$
(where $i$
is a natural number in  $[0,n-1]$).
The fundamental class of $Q_0$
modulo
$
W^{\langle k-1  \rangle}
$
equals to
$[p]$ by definition.
For $i>0$
we have:
\bqq
Q_i
=
T(\Sigma_i(p), -v)\cap F^{-1}([\l-i-1+a_k, \l-i])
\end{equation*}
Thus the homology class of $Q_i$
modulo its boundary equals to
the element
$\tau(\s_i(p))$
and
the application of \rref{f:sigma_r}
finishes the proof. $\qs$

Now we can define a map of $\wh L_-$-modules
\bq
C_*^-(f,v)\rTo^\xi \EE_*\tens{  L_-}\wh  L_-
\end{equation}
setting for the generators $p\in S_k(  F )\cap   W$
\bq
\label{f:xi_def}
\xi(p)
=
[p]
-
\sum_{r=0}^\infty \wi\tau(\HH^r_{k-1}(\s_1(p)))
\end{equation}
By \ref{c:heeq}
there is a chain \heq~
\bq
\l: \EE_*\tens{  L_-}\wh  L_-
\to
C_*^s(  V^-)\tens{  L_-}\wh  L_-
\end{equation}
By definition (\cite{pafest})
the chain homotopy equivalence
$\psi$
is
the composition
$\l\circ\xi$.
In \cite{pafest}
the fact that $\xi$
is a chain map is verified by a direct computation.
Using the methods developed in this paper
we shall  give a simpler proof, which constitutes the first part
of the Proposition
\ref{p:ksi_psi}
in the next section. 
The second part of this
proposition implies Theorem B.

\section{Proof of Theorem B}
\label{s:proofb}

In the first subsection we prove Theorem B
for $C^0$-generic gradients in $\GG(f)$.
The subsections \ref{su:inv_tors},
\ref{su:inv_zeta}
 contain the proof of the general case, which
 will be deduced from the $C^0$-generic case
by a perturbation argument. Here is the schema of this argument.
Let $v$ be an arbitrary Kupka-Smale gradient, and
$$
\tau(v)=w(M,f,v), \quad \z(v)=(\zeta_L(-v))^{-1}\in \ZZZ[[t]]
$$
be the corresponding invariants.
It suffices to prove that for every $n\geq 0$
the images
$\tau_n(v), \z_n(v)$
of these elements in the quotient ring
$L_n=\ZZZ[[t]]/t^n$
are equal.
Using Theorem \ref{t:cc}
pick a $C^0$-small perturbation $w$ of $v$
\sut~$w\in \GG_0(f)$.
We show that
$\tau_n(v)=\tau_n(w)$
(Section \ref{su:inv_tors} below),
and
$\z_n(v)=\z_n(w)$
(Section \ref{su:inv_zeta} below).
 Theorem B follows, since
$\z(w)=\tau(w)$
by the results of Subsection
\ref{su:c0_gener}.

\subsection{The case of $C^0$-generic $f$-gradient}
\label{su:c0_gener}

In this and the next  subsection
we work with the terminology of
Subsection \ref{su:psi}.
Thus $v$ is an $f$-gradient
satisfying the condition $(\gC)$,
$F:\bar M\to\RRR$ is a lift
of $f:M\to S^1$ and
$\l$ is a regular value of $F$.
We fix $\l$ up to the end of the section, and
the chain complex $C_*^-(f,v;\l)$
(see Subsection \ref{su:morse_wn}
for definition)
will be denoted $C_*^-(f,v)$ for brevity.
Recall the chain equivalence
$\mu:
C_*^-(f,v)
\to
 C_*^s( V^-)\tens{ L_-} \wh L_-
$
from Proposition \ref{p:mu}
and the homomorphism
\rref{f:xi_def}
of $\wh L$-modules
$\xi:C_*^-(f,v)\to \EE_*\tens{L_-}\wh L_-$.

\bepr\label{p:ksi_psi}
\been\item
$\xi$ is a chain map.
\item
The following diagram is chain homotopy commutative
\pa\pa
\begin{diagram}[size=2em,LaTeXeqno]
C_*^-(f,v)  &  &  & &
\\
\dTo_{\mu} &  & \rdTo(2,1)^{\xi}  & & \EE_*\tens{ L_-}\wh L_-
 \quad\lb{f:triangle}
\\
 C_*^s( V^-)\tens{ L_-} \wh L_-  & &   \ldTo(2,1)_{\l}      &   & \\
\end{diagram}
\pa\pa
\enen
\enpr
\Prf
Consider the map
\bq
\xi_n=\xi\otimes L_n:
C_*^-(f,v)
\Big/
t^n C_*^-(f,v)
\rTo\EE_*/t^n\EE_*
\end{equation}
Both the target and the source of $\xi_n$
are adjoint chain complexes corresponding to two different
filtrations of the pair
$
( V^-, t^n  V^-)$.
The first filtration is associated
 with a $t$-ordered Morse function
$\chi$ on $W_n$
(see Subsection \ref{su:morse_wn}),
 the second is the filtration
\bq
\Big\{ \Big( V^-_{\langle k  \rangle}
\cup t^n V^-
,
t^n  V^-
\Big)
\Big\}
\end{equation}
constructed in \ref{su:fil_cov}.
According to \ref{p:thin_fi}
 we can choose a $t$-ordered function $\chi:W_n\to\RRR$
and an ordered sequence
$(a_0,... ,a_{m+1})$
for $\chi$
in such a way that for every $k$ we have
\bq
\chi^{-1}([a_0, a_{k+1}])\cup t^n V^-
\sbs
V_{\langle k \rangle}^-\cup t^n V^-
\end{equation}
Thus the identity map of
$( V^-, t^n  V^-)$
sends the
terms
of the first filtration to the
respective  terms of
the second filtration.
It follows from the definition of
$\xi$
and
Proposition \ref{p:compu_d}
that the map
$\xi_n$
is exactly the map induced in the adjoint  complexes
by the identity map.
Thus $\xi_n$
is the chain map for every $n$, and therefore
 $\xi$ is also a chain map.
Moreover the diagram
obtained from
(\ref{f:triangle})
by tensoring with $ L_n$
is homotopy commutative by \ref{c:heeq}.
Therefore
(\ref{f:triangle})
is homotopy commutative by Proposition  \ref{p:noet}.
$\qs$

The case of $C^0$-generic gradients will be over with the next proposition.
\bepr\label{p:B_gener}
Let $v\in \GG(f)$. Then
\bq\label{f:B_gener}
\tau(v)=\z(v).
\end{equation}
\enpr
\Prf
In order to compute the torsion
of the chain map
$$
\Phi_\D=\chi^{-1}_\D\circ \Phi(M,f,v)
$$
it suffices to compute the torsion of the chain equivalence
$$
C_*^-(f,v)
\rTo^\xi \EE_*\tens{L_-}\wh L_-
$$
since the torsion of
$\chi^{-1}_\D\circ\l$
vanishes (see the remark just before 
Subsection \ref{su:psi}).
To this end note that for every generator
$\{p\}$
of
$C_*^-(f,v)$
corresponding to a critical point
$p\in S(f)$
its image 
$\xi(p)$
differs from the free generator
$[p]$
of $\wh\EE_*=\EE_*\tens{L_-}\wh L_-$
by an element in
$t\cdot\wh\EE_*$.
Therefore for every $k$ the family
$$
\{
[p]_{p\in S_k(\phi_1)}
,~
\xi([q])_{q\in S_k(\phi_1)}
,~
\wi\tau([r])_{r\in S_{k-1}(\phi_1)}
\}
$$
is a free base of 
$\wh\EE_k$, and the torsion of
$\xi$
is equal to the torsion of the acyclic complex
$$
\wh\EE_*'
=
\wh\EE_*\big/\Im\xi
$$
The module $\wh\EE_k'$
is isomorphic to
$
\wh\RR_k
\oplus\wh\RR_{k-1}$
where
$$
\wh\RR_k
=
\RR_k
\tens{L_-}\wh L_-
=
C_k(\phi_1, u_1)\tens{\ZZZ}\wh L_-.
$$
Therefore
the $\wh L_-$-base of 
$\wh\EE_*'
$
is identified with 
$S_k(\phi_1)\sqcup S_{k-1}(\phi_1)$
and the matrix of the boundary 
operator  
$$
\pr_k'
:
\wh\EE_k'
\to
\wh\EE_{k-1}'
$$
\wrt~this base 
is obtained immediately from \rref{f:bou_v}:
\bq\label{f:new_d}
\pr'_k=
\left(
\begin{matrix}
 \pr_{k}^{(1)}   & 1 -\HH_{k-1} \\
0 &  \d_{k-1} \\
\end{matrix}   \right)
\end{equation}
where
$$\d_{k}:
\wh\RR_k\to\wh\RR_{k-1}
$$
is some \ho~with $\d_{k-1}\circ \d_k =0$.
Recall that for erevy $k$ we have $\Im\HH_{k}\sbs t\cdot 
\RR_{k}
$
so that the homomorphism
$$
\big(1-\HH_{k}\big)
\otimes \id:
\RR_{k}\tens{L_-}\wh L_-
\to
\RR_{k}\tens{L_-}\wh L_-
$$
is invertible, and determines an element
$ [1-\HH_{k}]$
in $\ove{K}_1(\wh L_-)$.
The following lemma
is purely algebraic, its proof uses only the algebraic
information on the complex $\wh\EE'_*$, cited above.
\bele[\cite{pafest}, \S 6]\label{l:tors_E'}  
The torsion of the chain complex
$\wh\EE_*'$
is equal to
\bq
\sum_k
(-1)^{k}  [1-\HH_{k}].
\qquad \qs
\end{equation}
\enle
Now recall the formula 
\rref{f:zeta_and_h}
and the proposition is proved.
$\qs$

\subsection{Invariance of $\tau_n(v)$ under
small  $C^0$-perturbations of $v$}
\label{su:inv_tors}

Let $\phi:W_n\to[a,b]$
be any $t$-ordered Morse function
on $W_n$, \sut~$v$ is also a $\phi$-gradient.
Choose an ordering sequence
$(a_0,...,a_{m+1})$
for $\phi$.
For each $p$ choose a \nei~
$U_p$
of $p$ in $M$, so small that 
$\pi^{-1}(U_p)\cap W_n\sbs
\phi^{-1}(]a_k, a_{k+1}[)$
where $\pi:\bar M\to M$
is the infinite cyclic covering and
$k=\ind p$.
Let $U=\cup_pU_p$.
The aim of the present subsection is the following proposition.
\bepr\label{p:inv_tau}
There is $\e>0$
\sut~for any 
$w\in G(f)$
with
$||w-v||<\e$
\footnote 
{We recall that the symbol
$||\cdot||$ denotes the $C^0$-norm.}
and
$w|U=v|U$ 
we have
$$
\tau_n(w)=\tau_n(v).
$$
\enpr
\Prf
Let
$(X^{(i)}, t^n V^-)$ be the filtration of the pair
$(V^-, t^nV^-)$
corresponding to the function $\phi$
(see \rref{f:x_i}).
There is an isomorphism of the truncated Novikov complex
$ C_*^-(f,v) /   t^n    C_*^-(f,v) $
to
$C_*^{gr}(V^-, t^nV^-)$.
This isomorphism is base preserving, if we choose in
$C_*^{gr}(V^-, t^nV^-)$
the base formed by the homology classes of
the descending discs of the critical points.
Choose a triangulation $\D$
in such a way that $V$ is a simplicial subcomplex.
Then $W_n$ is a simplicial subcomplex of $\bar M$
and the  simplicial chain complex
$C_*^\D(V^-, t^nV^-)$
obtains a structure of a free based
$\ZZZ[t]/t^n$-module.
The element
$\tau_n(v)\in \ZZZ[t]/t^n$
(modulo $\pm 1$)
is the determinant of the torsion of the following composition:
 \bq
C_*^{gr}(V^-, t^nV^-)
\rTo^{J_n(v)}
C_*^s (V^-, t^nV^-)
\rTo
C_*^\triangle (V^-, t^nV^-)
\end{equation}
where the base in the left hand side
complex is formed by the homotopy classes of the descending discs
of $v$.
It is easy to prove that 
there is $\e>0$ \sut~every $f$-gradient $w$
with $||v-w||<\e$
and with
$w|U=v|U$ 
satisfies the two following conditions
\been\item
 $w$ is still a $\phi$-gradient.
\item
the homotopy classes of the descending discs corresponding to $v$
coincide with those corresponding to $w$.
\enen
Thus the chain equivalences $J(v)$ and $J(w)$
have the same domain and are chain homotopic.
Moreover, the bases in their domains are the same.
The proposition follows. $\qs$

\subsection{Invariance of $\z_n(v)$ under small $C^0$-perturbations of $v$}
\label{su:inv_zeta}

This part is more delicate.
Let $k\in\NNN$
and recall  the  (not everywhere defined) map
$(\stv)^k : V\to t^kV$.
The map
$t^{-k}(\stv)^k$
is a diffeomorphism
of an open subset in $V$ to another open subset of $V$.
For $v\in \GG(f)$
the set
$\FF_k(v)$
of the fixed points of this diffeomorphism
is finite.
Let $L_k(v)$
be the algebraic number of these fixed points, i.e.
$L_k(v)=\sum\varepsilon(p)$
where
$\varepsilon(p)=\pm 1$
stands for the index of the fixed point $p$.
Then we have the following formula (which follows
 easily from the definition):
\bq\label{f:delta}
\eta_L(-v) = \sum_{k\geq 1} \frac {L_k(v)}{k} t^k
\end{equation}
The next proposition asserts the invariance of $\zeta_n(v)$
\wrt~small
$C^0$-perturbations.
\bepr
\label{p:inv_zeta}
Let $v\in\GG(f)$, and $k\in\NNN$.
There is $\e>0$
\sut~ for every
$w\in \GG(f)$
with $||w-v||<\e$
we have
\bq\label{f:inv_zeta}
L_k(v)=L_k(w).
\end{equation}
\enpr
The proof occupies the rest of this section.
We shall work in the cobordism $W_k$.
Pick a number 
 $\d$
 sufficiently small so that
 $v~|~W_k$ is $\d$-separated.
Let
 \bqq
K^+(\d)=
\bigcup_{p\in S(F)\cap W_k} D_\d(p,-v),
\qquad
K^-(\d)=
\bigcup_{p\in S(F)\cap W_k} D_\d(p,v)
\end{equation*}

\bepr\label{p:disj_v}
For any sufficiently small $\d>0$
we have:
\bq\label{f:disj_v}
\FF_k(v)\cap K^+(\d)=\emptyset.
\end{equation}
\enpr
\Prf
Let
$K_r^+$
be the union of all ascending discs of $v$ corresponding to the critical points of $v$
of indices $\geq r$, i.e.
\bqq
K_r^+=
\bigcup_{ \substack{ p\in S(F)\cap W_k, \\ \ind p\geq r}} D(p,-v)
\end{equation*}
The set
$K_r^+$
is compact, as well as its $\d$-thickening $K_r^+(\d)
$ introduced in the next formula
\bq\label{f:kr_dplus}
K_r^+(\d)=
\bigcup_{\substack{ p\in S(F)\cap W_k,\\
 \ind p\geq r}} D_\d(p,-v),
\end{equation}
Similary, set
\bq\label{f:kr_min}
K_r^-=
\bigcup_{\substack{ p\in S(F)\cap W_k, \\  \ind p\leq r}} D(p,v),
\qquad
K_r^-(\d)=
\bigcup_{  \substack{p\in S(F)\cap W_k, \\ \ind p\leq r}} D_\d(p,v)
\end{equation}

The transversality condition implies
\bq
t^k(K_r^+)
\cap
K_r^-
=\emptyset.
\end{equation}
Since
$\kpd$, resp $\kmd$
form a fundamental system of \nei s
of
$\kpr$, resp. $\kmr$,
(see Proposition \ref{p:ddd})
there is $\d>0$
\sut~
\bq\label{f:***}
t^k(\kpd)\cap \kmd =\emptyset.
\end{equation}
Take $\d$ so small that
(\ref{f:***})
holds for all $r$.
Note that with this $\d$
we have
$\FF_k(v)\cap K^+(\d) = \vide$.
Indeed, let
$x\in \FF_k(v)$
and assume $x\in D_\d(p, -v)$
with $\ind p=r$.
Then
$x\in \kpd$
and thus $t^k x\notin \kmd$.
On the other hand the $(-v)$-trajectory
$\g(x,\cdot;-v)$
reaches $t^kV$
and intersects it
at the point 
$$
t^k x\in
 D_\d(p,v)
\cap
t^kV
\sbs
K_r^-(\d),
$$
contradiction. $\qs$

Our next aim is to show that the formula
\rref{f:disj_v}
still holds if we replace $v$ by its $C^0$-small perturbation, i.e. an $f$-gradient
$w$, which is close to $v$ in $C^0$-topology.
Choose $\d>0$
\sut~
the conclusion of Proposition \ref{p:disj_v}
hold.
Let $\d>\eta>\mu>0$.
Let $w$ be an $f$-gradient.
Similarly to the definitions above set
\bq\label{f:kr_w}
K^+(\d,w)=
\bigcup_{ p\in S(F)\cap W_k} D_\d(p,-w),
\qquad
K^-(\d,w)=
\bigcup_{ p\in S(F)\cap W_k} D_\d(p,w)
\end{equation}

\bepr\label{p:disj_w}
There is $\e>0$, \sut~ for every
$w\in \GG(f)$ with $||w-v||<\e$
we have:
\bq\label{f:disj_w}
\FF_k(w)\cap K^+(\eta,w)=\vide.
\end{equation}
\enpr
\Prf
As in the proof of Proposition \ref{p:disj_v}
it suffices to show that for any given $r\in\NNN, 1\leq r\leq n$
and a $C^0$-small  perturbation $w$ of $v$ we have:
$$
t^k(K_r^+(\eta,w))
\cap
K_r^-(\eta,w)
=\vide.
$$
Here the sets $
K_r^+(\eta,w), ~K_r^-(\eta,w)$
are defined similarly to \rref{f:kr_dplus}, \rref{f:kr_min},
namely
\bqq
K_r^+(\eta,w)=
\bigcup_{\substack{ p\in S(F)\cap W_k, \\  \ind p\geq r}} D_\eta(p,-w),
\qquad
K_r^-(\eta,w)=
\bigcup_{\substack{ p\in S(F)\cap W_k, \\ \ind p\leq r}} D_\eta(p,w).
\end{equation*}
These sets 
are the $\eta$-thickenings of the following
unions of descending, respectively ascending discs:
\bqq
K_r^+(w)=
\bigcup_{  \substack{ p\in S(F)\cap W_k, \\ \ind p\geq r}} D(p,-w),
\qquad
K_r^-(w)=
\bigcup_{\substack{p\in S(F)\cap W_k, \\ \ind p\leq r}} D(p,w),
\end{equation*}
so that we have in particular
\begin{equation*}
K_r^+(w)
\sbs
\Int ( K_r^+(\eta,w) ) , \quad
K_r^-(w)
\sbs
\Int (   K_r^-(\eta,w)    ).
\end{equation*}
By Lemma \ref{l:cont_2} we have
\bqq
K_r^+(\eta,w)  \sbs   \Int (K_r^+(\d,v)),\qquad
K_r^-(\eta,w)  \sbs   \Int (K_r^-(\d,v)).
\end{equation*}
for every sufficiently $C^0$-small  perturbation
$w$ of $v$. Now apply (\ref{f:***})
and the proof is over. $\qs$

Next we shall interpret the number
$L_k(w)$
for $C^0$-small perturbations $w$ of the vector field $v$
in terms of
fixed point indices  of some continuous maps.
For this we need some preliminaries.
The next lemma follows from Lemma \ref{l:cont_2}.

\bele\label{l:v_w}
There is $\e>0$, \sut~ for every
$w\in G(f)$ with $||w-v||<\e$  we have
\bq\label{f:v_w}
K_r^+(w)\sbs \Int K_r^+(\mu,v). \qquad \square
\end{equation}
\enle

Pick any $(m-1)$-dimensional
 compact \sma~ $B$ with boundary of $V$, \sut~
\bq
K^+(\mu, v)\sbs \Int B \sbs B \sbs \Int K^+(\eta,v).
\end{equation}
\bele\label{l:l_B}
There is $\e>0$, \sut~ for every
$w\in G(f)$ with $||w-v||<\e$  we have
\bq\label{f:l_B}
K^+( w)\sbs \Int B \sbs B \sbs \Int K^+(\eta,w)
\end{equation}
\enle
\Prf
The first inclusion follows from Lemma \ref{l:v_w}.
The last one  follows from Lemma \ref{l:cont_1}.$\qs$

The set
$$
C= V\sm \Int B
$$
is a compact $(m-1)$-dimensional \sma~ with boundary
of $V$.
It follows from (\ref{f:l_B})
that every $(-w)$-trajectory starting at a point of $C$
reaches $t^kV$.
\bepr\label{p:w_C}
There is $\e>0$, \sut~ for every
$w\in G(f)$ with $||w-v||<\e$  we have
$$
\FF_k(w)\sbs \Int C.$$
\enpr
\Prf
Indeed
$\FF_k(w)\cap K^+(\eta,w)=\vide$
by Proposition \ref{p:disj_w}
therefore $\FF_k(w)\sbs V\sm  B$.$\qs$

Now let $w$ be an $f$-gradient \sut
~\rref{f:l_B}
holds.
The set $C$ is then 
in the domain of definition of the map
$(\stw)^k$
 and the restriction to $C$ of the map
$t^{-k}
(\stw)^k$
is a $\smo$ embedding $S_w:C\to V$.
Assuming further that $w\in \GG(f)$
and that \rref{f:disj_w} holds,
we deduce that the set
$\Fix(S_w)=\FF_k(w)$
is finite and does not intersect
with $\pr C$.
Denote the sum of the indices of all the fixed points
of $S_w$
by $\JJ(S_w)$,
then we have by definition
\bq\lb{f:equal}
L_k(w)=\JJ(S_w).
\end{equation}
Thus for such $f$-gradients $w$
the computation of the number
$L_k(w)$
is reduced to the computation of fixed point indices
of the continuous map
$S_w:C\to V$
associated to $w$.
Actually the map
$S_w:C\to V$
is defined for every $C^1$ vector field $w$
sufficiently close to $v$ in $C^0$ topology
and the correspondance 
$w\mapsto S_w$
is continuous (\wrt~ the $C^0$-topology on the set of
vector fields and on the set
of maps $C\to V$).
This is the contents of the next proposition, which follows directly from
\cite{paadv}, Prop. 2.66.
 \bepr
Let $K\sbs \pr_1 W$
be a compact set, and $u$ be a $C^1$ vector field on $W$, \sut~
$u|\pr_1 W$ points outward $W$, and $u|\pr_0 W$ points inward $W$.
Assume that every $(-u)$-trajectory starting at $K$
reaches $\dow$, and denote the restriction of
$\stu$ to $K$
by $S_u$.
Then there is $\e>0$
\sut~for every
$C^1$ vector field $\wi u$ on $W$ with
$||\wi u -u ||<\e$, every
$(-\wi u)$-trajectory
starting at $K$ reaches $\dow$
and the map
$\wi u\mapsto S(\wi u)$
is continuous. $\qs$
\enpr

A continuous map
$g:C\to V$
will be called {\it of class $\FF$}
if the fixed point set $\Fix(g)$
is finite and $\Fix(g)\cap \pr C=\emptyset$.
For a map
$g$
of class $\FF$
we denote the sum of indices of its fixed points by
$\JJ(g)$.
Our main example of maps of class $\FF$
is any map $S_w$
with $w\in \GG(f)$ and $||w-v||<\e$,
where $\e>0$
is chosen so small that 
Proposition \ref{p:disj_w}, and Lemmas
\ref{l:v_w}, \ref{f:l_B}
hold with this $\e$.
Now the proposition \ref{p:inv_zeta}
follows from \rref{f:equal} and the next lemma.
\bele\label{l:cont_ind}
Let $g:C\to V$
be a continuous map of class $\FF$.
There is
$\d>0$
\sut~for every map $h:C\to V$
of class $\FF$
with
$d(h,g)<\d$
we have
$$
\JJ(h)
=
\JJ(g).
$$
\enle
\Prf
Let
$p\in \Fix(g)$.
Choose a chart
$\Psi:U\to\RRR^n$
of the \ma~$V$ 
with $p\in U$ and $\Psi(p)=0$.
The inverse image of the closed
disc
$D(0,R)\sbs \RRR^n$
\wrt~$\Psi$
will be denoted $D_R(p)$.
Similarly
the inverse image of the open
disc
$B(0,R)\sbs \RRR^n$
\wrt~$\Psi$
will be denoted $B_R(p)$.
Choose $R$ so small that for every $p$ the only fixed point
of $g$ in $D_R(p)$
is $p$, and the discs $D_R(p)$
are disjoint for different points $p$.
Pick a number $\l\in ]0, R[$
so small that
for each $p\in \Fix(g)$
we have
$g(D_\l(p))\sbs B_R(p)$.
The next lemma is proved by a 
standard continuity argument.

\bele\lb{l:s_nu}
There is $\nu>0$
\sut~for every map
$h:C\to V$
of class $\FF$, \sut~$d(h,g)<\nu$
we have:
\bq\lb{f:s_nu}
\left\{
\begin{aligned}   
h(D_\l(p))    &\sbs B_R(p) \quad \mx{ for every } p\in \Fix(g) \\
\Fix(h)      &\sbs \cup_p D_{\l/2}(p). \hspace*{5cm}\qs
\end{aligned}
\right.
\end{equation}
\enle
For every $h$ satisfying \rref{f:s_nu}
the set $\Fix(h)$
falls into the disjoint union 
of the sets of the fixed points of
the maps 
$$
h[p]=h~|~D_\l(p)
:
D_\l(p)
\to
B_R(p).
$$

Thus it remains  to prove that 
there is $\rho$ \sut~for every $h$ with 
$d(h,g)<\rho$ 
we have 
\bq\label{f:equa_ind}
\JJ(h[p])
=
\JJ(g[p])\quad 
\mx{~for every }~ p
\end{equation}
Let $\mu$
be the injectivity radius of the riemannian 
manifold $V$, so that
for every $\mu'<\mu$
every two maps
$F_1, F_2:C\to V$
satisfying
\bq\label{f:dist}
d(F_1, F_2)\leq \mu'
\end{equation}
are homotopic in the class of maps satisfying
\rref{f:dist}.
Then for every
$h:C\to V$
satisfying
$$
d(h,g)\leq \min(\mu/2, \nu)
$$
the maps
$g[p], h[p]:D_\l(p)\to   B_R(p)$
are homotopic
via a homotopy
$H_t:D_\l(p)\to    B_R(p)$,
\sut~for every
$t$ the set
of fixed points of $H_t$
is compact and
contained in $D_{\l/2}(p)$.
Then
the equality 
\rref{f:equa_ind}
holds 
by the homotopy invariance of the fixed point index
(\cite{dold}, 5.8), and
the proof of the lemma is over. $ \qs$

Now we can finish the proof of Theorem B.
As we have mentioned already it suffices to prove that for every
$n\in\NNN$
we have
$\z_n(v)=\tau_n(v)$.
For each $p\in S(f)$
choose a \nei~$U_p$
of $p$
as in the beginning of Subsection \ref{su:inv_tors}.
Choose $\e>0$
so small that for any
$w\in G(f)$
with
$||w-v||<\e$
the conclusions of propositions \ref{p:inv_tau},
\ref{p:inv_zeta}
hold.
Proposition \ref{t:cc}
implies the existence of an $f$-gradient
$w\in GC(f)$
\sut
~$||v-w||<\e/3$.
By a standard Kupka-Smale argument
there is an $f$-gradient
$u\in \GG(f)$
with
$||u-w||<\e/3$.
Moreover the gradient $u$
is still in $GC(f)$
if only 
~$||u-w||$
is sufficiently small
(since $GC(f)$
is $C^0$-open).
Proposition \ref{p:B_gener}
implies
$\tau(u)=\zeta(u)$;
in particular we have
$\tau_n(u)=\z_n(u)$
for every $n$.
Note finally that Proposition
\ref{p:inv_zeta}
implies
$\z_n(u)=\z_n(v)$
and Proposition
\ref{p:inv_tau}
implies
$\tau_n(u)=\tau_n(v)$.
The proof of Theorem B is over. $\qs$

\newpage

\label{figure}

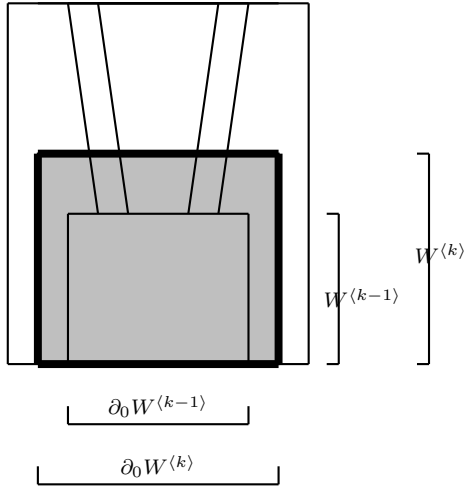
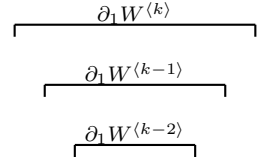
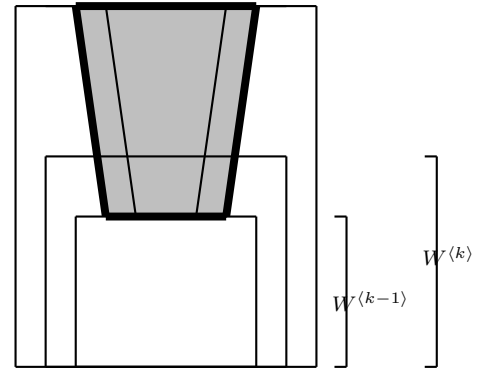
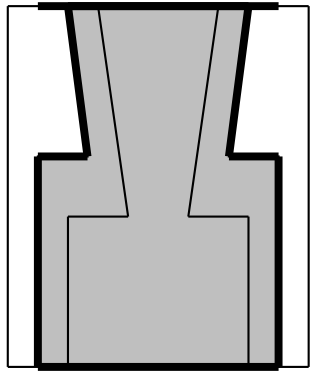
\psset{xunit=.4cm,yunit=.4cm}
\psset{runit=.4cm}

%\psgrid(0,0)(0,0)(12,12)

\begin{pspicture}(0,0)(20,20)

%\input  zk.doc

%bol'shoi kvadr
\pscustom[fillstyle=solid, fillcolor=lightgray]
{%
\psline(1,0)(1,7)
\psline(1,7)(9,7)
\psline(9,7)(9,0)
\psline(9,0)(1,0)
}

\psline[linewidth=3pt](1,0)(1,7)
\psline[linewidth=3pt](1,7)(9,7)
\psline[linewidth=3pt](9,7)(9,0)
\psline[linewidth=3pt](9,0)(1,0)

%malyi kvadr
 \psline(2,0)(2,5)
\psline(2,5)(8,5)
\psline(8,5)(8,0)
\psline(8,0)(2,0)

%bol'shaya trap.
\psline(3,5)(2,12)
\psline(2,12)(8,12)
\psline(8,12)(7,5)

%malaya trap.
\psline(4,5)(3,12)
\psline(3,12)(7,12)
\psline(7,12)(6,5)

%verh krysha
\psline(1,12)(9,12)

%obschaya ramka
\psline(0,0)(0,12)
\psline(0,12)(10,12)
\psline(10,12)(10,0)
\psline(10,0)(0,0)

%skobki sboku
\psline(10.6,0)(11,0)
\psline(11,0)(11,5)
\psline(11,5)(10.6,5)
%vtoraya
\psline(13.6,0)(14,0)
\psline(14,0)(14,7)
\psline(14,7)(13.6,7)

%bykvy sboku
\rput(11.8,2.2){$ {\scriptstyle W^{\langle k-1\rangle}}$}
\rput(14.4,3.7){${\scriptstyle W^{\langle k\rangle}}$}

%nizhnie  nadpisi
\rput(5,-1.4){${\scriptstyle \pr_0 W^{\langle k-1\rangle}}$}
\rput(5,-3.4){${\scriptstyle 
\pr_0 W^{\langle k\rangle}}$}

%skobka snizu pervaya
\psline(2,-2)(2,-1.4)
\psline(2,-2)(8,-2)
\psline(8,-2)(8,-1.4)

%skobka snizu vtoraya
\psline(1,-4)(1,-3.4)
\psline(1,-4)(9,-4)
\psline(9,-4)(9,-3.4)

%aixs_y
\psline(16,-1)(16,13)
\psdots[dotsize=4pt](16,0)(16,5)(16,7)(16,12)
\rput(16.7,0){$ {\scriptstyle a}$}
\rput(16.7,5){${\scriptstyle a_k}$}
\rput(16.9,7){${\scriptstyle a_{k+1}}$}
\rput(16.7,12){${  \scriptstyle  b}$}

\rput(5,15)
{ 
{The set ${Z_k}$  
(shaded).}
}

\end{pspicture}

\begin{pspicture}(-21,-1.3)(20,-1.3)

%\input  tk.doc

%bol'shaya trap.
\pscustom[fillstyle=solid, fillcolor=lightgray]
{%
\psline(3,5)(2,12)
\psline(2,12)(8,12)
\psline(8,12)(7,5)
\psline(7,5)(3,5)
}

\psline[linewidth=3pt](3,5)(2,12)
\psline[linewidth=3pt](2,12)(8,12)
\psline[linewidth=3pt](8,12)(7,5)
\psline[linewidth=3pt](7,5)(3,5)

%bol'shoi kvadr
\psline(1,0)(1,7)
\psline(1,7)(9,7)
\psline(9,7)(9,0)
\psline(9,0)(1,0)

%malyi kvadr
\psline(2,0)(2,5)
\psline(2,5)(8,5)
\psline(8,5)(8,0)
\psline(8,0)(2,0)

%malaya trap.
\psline(4,5)(3,12)
\psline(3,12)(7,12)
\psline(7,12)(6,5)

%verh krysha
\psline(1,12)(9,12)

%obschaya ramka
\psline(0,0)(0,12)
\psline(0,12)(10,12)
\psline(10,12)(10,0)
\psline(10,0)(0,0)

%skobki sboku
\psline(10.6,0)(11,0)
\psline(11,0)(11,5)
\psline(11,5)(10.6,5)
%vtoraya
\psline(13.6,0)(14,0)
\psline(14,0)(14,7)
\psline(14,7)(13.6,7)

%bykvy sboku
\rput(11.8,2.2){$ {\scriptstyle W^{\langle k-1\rangle}}$}
\rput(14.4,3.7){${\scriptstyle W^{\langle k\rangle}}$}

%skobki sverhu
\psline(3,14)(3,13.6)
\psline(3,14)(7,14)
\psline(7,14)(7,13.6)

%skobki sverhu vtoraya
\psline(2,16)(2,15.6)
\psline(2,16)(8,16)
\psline(8,16)(8,15.6)

%skobki sverhu tretia
\psline(1,18)(1,17.6)
\psline(1,18)(9,18)
\psline(9,18)(9,17.6)

%verhnie nadpisi
\rput(5,14.4){$ {\scriptstyle \pr_1 W^{\langle k-2\rangle}}$}
\rput(5,16.4){${\scriptstyle \pr_1 W^{\langle k-1\rangle}}$}
\rput(5,18.4){${\scriptstyle \pr_1 W^{\langle k\rangle}}$}

\rput(5,-2)
{ {The set ${T_k}$  {(shaded)}.}}

\end{pspicture}

\begin{pspicture}(0,0)(20,20)

%\input  wk.doc

%obschaya ramka
\psline(0,0)(0,12)
\psline(0,12)(10,12)
\psline(10,12)(10,0)
\psline(10,0)(0,0)

%aixs_y
\psline(16,-1)(16,13)
\psdots[dotsize=4pt](16,0)(16,5)(16,7)(16,12)
\rput(16.7,0){$ {\scriptstyle a}$}
\rput(16.7,5){${\scriptstyle a_k}$}
\rput(16.9,7){${\scriptstyle a_{k+1}}$}
\rput(16.7,12){${  \scriptstyle  b}$}

\rput(5,-2)
{ The set ${
 W^{\langle k\rangle}
}$ 
  is shaded
}

\rput(5,-4)
{
and the set 
${
 W^{\langle k-1\rangle}
}$ 
is shown by thin lines.
}

\rput(15,-10){
Handle-like filtration of the cobordism $W$
(Subsection \ref{su:c_hand}).
}

\pscustom[fillstyle=solid, fillcolor=lightgray]
{%
\psline(1,0)(1,7)
\psline(1,7)(2.65,7)
\psline(2.65,7)(2,12)
\psline(2,12)(8,12)
\psline(8,12)(7.35,7)
\psline(7.35,7)(9,7)
\psline(9,7)(9,0)
\psline(9,0)(1,0)
}

\psline[linewidth=3pt](1,0)(1,7)
\psline[linewidth=3pt](1,7)(2.65,7)
\psline[linewidth=3pt](2.65,7)(2,12)
\psline[linewidth=3pt](2,12)(8,12)
\psline[linewidth=3pt](8,12)(7.35,7)
\psline[linewidth=3pt](7.35,7)(9,7)
\psline[linewidth=3pt](9,7)(9,0)
\psline[linewidth=3pt](9,0)(1,0)
\psline[linewidth=3pt](1,12)(9,12)

\psline(2,0)(2,5)
\psline(2,5)(4,5)
\psline(4,5)(3,12)
\psline(3,12)(7,12)
\psline(7,12)(6,5)
\psline(6,5)(8,5)
\psline(8,5)(8,0)
\psline(8,0)(2,0)

\end{pspicture}

\begin{pspicture}(-21,-1)(20,-1)

%\input  wkm.doc

%malyi kvadr
 \psline(2,0)(2,5)
\psline(2,5)(8,5)
\psline(8,5)(8,0)
\psline(8,0)(2,0)

%obschaya ramka
\psline(0,0)(0,12)
\psline(0,12)(10,12)
\psline(10,12)(10,0)
\psline(10,0)(0,0)

\rput(5,-2)
{ {The set ${
W^{\langle k-1\rangle}
}$  {(shaded)}.}}

\pscustom[fillstyle=solid, fillcolor=lightgray]
{%
\psline(2,0)(2,5)
\psline(2,5)(4,5)
\psline(4,5)(3,12)
\psline(3,12)(7,12)
\psline(7,12)(6,5)
\psline(6,5)(8,5)
\psline(8,5)(8,0)
\psline(8,0)(2,0)
}

%[linewidth=3pt]

\psline[linewidth=3pt](2,0)(2,5)
\psline[linewidth=3pt](2,5)(4,5)
\psline[linewidth=3pt](4,5)(3,12)
\psline[linewidth=3pt](3,12)(7,12)
\psline[linewidth=3pt](7,12)(6,5)
\psline[linewidth=3pt](6,5)(8,5)
\psline[linewidth=3pt](8,5)(8,0)
\psline[linewidth=3pt](8,0)(2,0)
\psline[linewidth=3pt](2,12)(8,12)

\end{pspicture}

\newpage
\label{refer}

\end{document}